\documentclass[a4paper]{article}

\usepackage[margin=1in]{geometry} 

\usepackage[utf8]{inputenc}
\usepackage[T1]{fontenc}

\usepackage{subfigure}
\usepackage{mathtools}
\usepackage{amsmath}
\usepackage{amsthm}
\usepackage{amssymb}
\usepackage{physics}
\usepackage{float}
\usepackage[colorlinks=true]{hyperref}
\usepackage{booktabs}
\usepackage{array}
\usepackage{bm}

\usepackage{amsmath,amsthm,amssymb}
\usepackage[ruled,vlined]{algorithm2e}

\usepackage{threeparttable}
\usepackage{multirow}
\usepackage{booktabs}
\usepackage{color}
\usepackage{url}
\usepackage{lineno}

\hypersetup{%
    unicode,
    pdfauthor={Author One, Author Two, Author Three},
    pdftitle={TSGO},
    pdfkeywords={Continuous optimization},
}
\usepackage{authblk}
\usepackage{multirow}

\usepackage[sort&compress,numbers,square]{natbib}
\usepackage{graphicx, color}

\definecolor{OurRed}{RGB}{184,80,75}
\definecolor{OurBlue}{RGB}{88, 118, 185}
\definecolor{OurPurple}{RGB}{103, 82, 158}
\hypersetup{%
    colorlinks=true,
    linkcolor=OurBlue,
    citecolor=OurRed,
    urlcolor=OurPurple,
    pdftitle={An efficient optimization model and tabu search-based global optimization approach for continuous $p$-dispersion problem},
    pdfnewwindow=true,
}

\title{An efficient optimization model and tabu search-based global optimization approach for continuous $p$-dispersion problem}

\author{Xiangjing Lai}
\affil{Institute of Advanced Technology, Nanjing University of Posts and Telecommunications, Nanjing 210023, P.R. China; email:laixiangjing@gmail.com}
\author{Zhenheng Lin}
\affil{Institute of Advanced Technology, Nanjing University of Posts and Telecommunications, Nanjing 210023, P.R. China; email:linzhenheng@outlook.com}
\author{Jin-Kao Hao\thanks{jin-kao.hao@univ-angers.fr}}
\affil{LERIA, Universit$\acute{e}$ d'Angers, 2 Boulevard Lavoisier, 49045 Angers, France; email: jin-kao.hao@univ-angers.fr}
\author{Qinghua Wu}
\affil{School of Management, Huazhong University of Science and Technology, 430074 Wuhan, P.R.China; email:
qinghuawu1005@gmail.com}

\begin{document}
\maketitle

\begin{abstract}
Continuous $p$-dispersion problems with and without boundary constraints are NP-hard optimization problems with numerous real-world applications, notably in facility location and circle packing, which are widely studied in mathematics and operations research. In this work, we concentrate on general cases with a non-convex multiply-connected region that are rarely studied in the literature due to their intractability and the absence of an efficient optimization model. Using the penalty function approach, we design a unified and almost everywhere differentiable optimization model for these complex problems and propose a tabu search-based global optimization (TSGO) algorithm for solving them. Computational results over a variety of benchmark instances show that the proposed model works very well, allowing popular local optimization methods (e.g., the quasi-Newton methods and the conjugate gradient methods) to reach high-precision solutions due to the differentiability of the model. These results further demonstrate that the proposed TSGO algorithm is very efficient and significantly outperforms several popular global optimization algorithms in the literature, improving the best-known solutions for several existing instances in a short computational time. Experimental analyses are conducted to show the influence of several key ingredients of the algorithm on computational performance.

\end{abstract}

\section{Introduction}
\label{intro}

Equal circle packing and point arrangement are two important and closely related global optimization problems in mathematics and operations research. The equal circle packing problem consists of packing a fixed number $p$ of non-overlapping congruent circles into a bounded region to maximize the common radius of the circles \citep{addis2008disk,lopez2011heuristic,melissen1993densest,szabo2007new,weaire2008pursuit}. As a general model, the problem has a variety of practical applications ranging from circular cutting to container loading \citep{castillo2008solving}. By contrast, the point arrangement problem aims to scatter $p$ given points in a bounded region to maxmize the minimum distance between the points \citep{akiyama2002maximin,costa2013valid,goldberg1970packing,melissen1993densest,schwartz1970separating}. This problem finds application in areas such as facility location and design of computer experiments \citep{dimnaku2005approximate,DREZNER2019105}. For example, in the facility location problem \citep{DREZNER2019105,drezner1995solving,lopez2021new}, a facility is represented by a dispersion point and the goal is to place $p$ obnoxious facilities (e.g., nuclear facilities) in locations to maximize the minimum distance between them. These two problems can be mutually converted into each other if the given bounded region is a regular convex set (e.g., a convex polygon, a circle or the surface $S^2$ of a sphere), and should be treated separately otherwise.

The equal circle packing and point arrangement problems are also known under the general name of the continuous $p$-dispersion problem (CpDP) \citep{drezner1995solving,baur2001approximation,dimnaku2005approximate,chen2019efficient,dai2021repulsion} and can be described as follows. Given a positive integer $p$, a bounded region $\Omega$ composed of a container $C$ containing $K$ ($K \ge 0$) holes $\{U_1, U_2, \dots, U_K\}$ (i.e., $\Omega = C \setminus \cup_{k=1}^{K} U_{k}$), the CpDP is to place $p$ given dispersion points  $\{c_1,c_2,\dots, c_p\}$ in $\Omega$ so that the minimum distance from a dispersion point to the boundaries of the container and holes is no less than $D_b$, and the minimum distance $D$ between the dispersion points is maximized, where $D_b$ can be a non-negative constant or an increasing function with respect to $D$. Formally, the CpDP can be expressed as:

\begin{equation}\label{DMax}
\mathrm{Maximize} \quad D
\end{equation}
\begin{equation}\label{betweenpoints}
\mathrm{Subject \ to} \quad d(c_i,c_j) \ge D,  1 \le i \neq j \le p
\end{equation}
\begin{equation}\label{BoundaryC}
\quad  d(c_i,\partial C) \ge D_b, i = 1,2,\dots, p
\end{equation}
\begin{equation}\label{BoundaryH}
\quad  d(c_i,\partial U_{k}) \ge D_b, i = 1,2,\dots, p,  k = 1,2,\dots, K
\end{equation}
\begin{equation}\label{Containment}
\quad  c_i \in C \setminus \cup_{k=1}^{K} U_{k}, i = 1,2,\dots, p
\end{equation}
where $D$ is the objective function value to be maximized and represents the minimum distance between points, $d(c_i,c_j)$ is the Euclidean distance between points $c_i$ and $c_j$, $\partial C$ and $\partial U_{k}$ denote respectively the boundaries of $C$ and $U_{k}$, $d(c_i,\partial C)$ is the Euclidean distance between $c_i$ and $\partial C$ defined as $min \{d(c_i,c) : c \in \partial C\}$. $D_b$ denotes the allowed minimum distance between a dispersion point and the boundaries of the container and holes, and $D_b$ can be a constant or an increasing function with respect to $D$, such as $D_b(D) = \lambda \times D$  with $\lambda\in [0,0.5]$. The constraints (\ref{betweenpoints}) are called the separation constraints, the constraints (\ref{BoundaryC}) and (\ref{BoundaryH}) are called the boundary constraints, and the constraints (\ref{Containment}) are called the containment constraints. Following \citet{dai2021repulsion}, the container and holes are represented by a polygon because a simply-connected region can be approximated by a convex or non-convex polygon with multiple edges. The CpDP is a difficult nonlinear non-convex optimization problem like those in \citep{bierlaire2010heuristic,geissler2018solving,ugray2007scatter} and is computationally challenging especially for large-scale instances.

As mentioned above, depending on the setting of $D_b$ the CpDP problem formulated by Eqs. (\ref{DMax}--\ref{Containment}) has two well-known variants (i.e., the CpDP problems with and without a boundary constraint), corresponding to the equal circle packing problem and the point arrangement problem, respectively. First, if $D_b$ is $\frac{1}{2} \times D$, the corresponding problem is called the CpDP problem with a boundary constraint, where the dispersion points must be at least $D_b$ ($=\frac{D}{2}$) away from the boundaries of the container and  holes. Second, if $D_b$ is set to 0, the corresponding problem is called the CpDP problem without a boundary constraint, where the dispersion points are allowed to approach the boundaries of the container and holes. One observes that the only difference between the equal circle packing and point arrangement problems lies in whether the dispersion points (or the centers of circles) are allowed to approach the boundaries of the container and holes.

Due to their NP-hard feature \citep{demaine2010circle} and potential applications, the equal circle packing and point assignment problems have been widely studied in the literature by researchers from mathematics, facility location theory, and the design of computer experiments. From the computational and theoretical perspectives, there are a host of studies in the literature (e.g., see (\cite{Lai2022,szabo2007new}) for equal circle packing). We summarize the most representative studies as follows.

\citet{melissen1993densest} studied the equal circle packing problem and the point assignment problem in an equilateral triangle, while \citet{GRAHAM1998139} tackled the former problem in a circular region using a billiard simulation method. \citet{birgin2008minimizing} proposed a twice-differentiable nonlinear optimization model for the circle and sphere packing problems in a simple convex container. From contrasting perspectives, \citet{grosso2010solving} and \citet{addis2008disk} respectively used a monotonic basin hopping (MBH) algorithm and a population-based basin-hopping (PBH) algorithm for the problem of packing equal circles into a square and circular container. \citet{mladenovic2005reformulation} and \citet{lopez2011heuristic} designed heuristic algorithms based on the formulation space search method for the equal circle packing problem in a simple convex container. \citet{huang2010greedy} presented a greedy vacancy search algorithm for packing equal circles into a square container, while more recently, \citet{stoyan2020optimized} designed several heuristic strategies for the circle and sphere packing problems in a simple convex container. Very recently, \citet{amore2023circle} studied the equal circle packing problem in regular polygons via optimizing a related energy function with a heuristic algorithm. Complementing these studies, the popular Packomania website \citep{Specht2023} provides the best-known solutions collected from different researchers for the equal circle packing problems in a variety of simple convex containers. In the design of computer experiments, the point assignment problem was studied as the max-min distance design problem \citep{Mu2017,van2007maximin}.

Apart from the studies involving simple convex containers, several studies tackle the equal circle packing problem in a non-convex simply-connected region. For example, by the movement of circles in a rotated container under shaking and gravity, \citet{machchhar2017dense} proposed a two-phase heuristic algorithm for finding the densest packing of congruent circles in an arbitrary non-convex container with B-spline boundaries. Via the discretization of solution space and the popular bottom-left placement strategy, \citet{yuan2018packing} proposed an angle bisector heuristic algorithm for packing equal circles into a convex or concave polygonal container, reporting computational experiments that show their proposed algorithm can produce high-quality solutions for some small-scale instances. Recently, \citet{dai2021repulsion} introduced a heuristic algorithm for the continuous $p$-dispersion problem in a simply-connected non-convex polygonal container. The idea of their algorithm is to model the circles (or points) as physical bodies and induce the circles to move over time under repulsive forces between circles and between the circles and container boundaries. Computational results show that their algorithm significantly outperforms several previous algorithms including the angle bisector heuristic algorithm \citep{yuan2018packing}, establishing this algorithm as the current state-of-the-art algorithm for the continuous $p$-dispersion problem in a non-convex polygonal container.

Our literature review reveals that a large number of previous studies for the continuous $p$-dispersion problem focus on the case where the region is simply connected (in particular, involving a regular convex region such as a circular or square container). However, many practical situations involve a non-convex multiply-connected region, and the literature does not provide suitable optimization models nor efficient global optimization algorithms for dealing with the corresponding continuous $p$-dispersion problems, a limitation that provides the main motivation of this study. It should be clear that the non-convexity and multiply-connected feature of the given region makes the problem much more difficult to handle compared to the popular models involving a regular convex region.

The main contributions of our work can be summarized as follows:
\begin{enumerate}
  \item We extend the research of the continuous $p$-dispersion problem from a simply-connected region to non-convex multiply-connected regions in contrast to the focus of all previous studies on the simply-connected region case.
  \item We introduce a unified optimization model for the continuous $p$-dispersion problems with and without a boundary constraint, providing a model that can take advantage of popular continuous solvers such as the limit memory BFGS method to perform local optimization, thus making it possible for the global algorithm to reach high-precision solutions.
  \item We propose an effective global optimization algorithm called TSGO for the continuous $p$-dispersion problem both with and without a boundary constraint. The proposed algorithm is capable of addressing a variety of cases, including those with a convex, non-convex, simply-connected or a multiply-connected region. The source codes of the proposed algorithm are made publicly available for potential real-world applications.
  \item We provide a set of benchmark instances with a non-convex multiply-connected region, which fill the gap in the literature devoted to such instances. We expect that these instances will be useful for the evaluation and comparison of algorithms in the future.
\end{enumerate}

The rest of the paper is organized as follows. Section \ref{Models} presents a unified optimization model for the continuous $p$-dispersion problems in a non-convex multiply-connected region. Section \ref{Method} describes   the proposed TSGO algorithm. Section \ref{Experimental_results} evaluates the performance of the TSGO algorithm. Section \ref{discussion} presents an analysis to get insights into the performance of the algorithm. The last section summarizes the contributions of the present work and provides some perspectives for future studies. 

\section{A unified optimization model for the continuous $p$-dispersion problems}
\label{Models}

The goal of this section is to present a unified optimization model for the continuous $p$-dispersion problems with and without boundary constraints.

\subsection{Basic idea underlying the optimization model}
\label{IdeaOfModel}

The continuous $p$-dispersion problems studied are constrained max-min optimization problems, which are usually difficult to handle using popular local optimization methods (e.g., the quasi-Newton method). We adopt a popular approach to solving circle packing problems \citep{huang2010greedy,Lai2022,lai2023perturbation,lai2023iterated}, whose central idea is to transform the max-min optimization problem into a series of unconstrained optimization subproblems in which the allowed minimum distance between points is fixed to a constant $D$, and then solve them in sequence by an unconstrained global optimization algorithm.

The differentiability of the objective function is a very important feature, which allows a global optimization algorithm to reach high-precision solutions by using popular gradient-based local optimization methods, such as the quasi-Newton and conjugate gradient methods. However, for the studied continuous $p$-dispersion problem, it is very challenging to create a differentiable optimization model due to the non-convexity and multiply-connected feature of the region $\Omega$.

In view of the foregoing considerations, we use the idea of repulsions and attractions between the dispersion points and the boundaries of the container and holes to devise a unified unconstrained optimization model with an almost everywhere differentiable objective function for the continuous $p$-dispersion problems with and without boundary constraints, where the allowed minimum distance between dispersion points is fixed to a constant $D$. The following subsections describe this optimization model in detail.

\subsection{Continuous $p$-dispersion problem with an allowed minimum distance $D$ between points}
\label{Model1}

The continuous $p$-dispersion problem with an allowed minimum distance $D$ between points (CpDP-D for short) is a subproblem of the continuous $p$-dispersion problem studied in this work, whose goal is to find a feasible solution satisfying the constraints described by Eqs. (\ref{betweenpoints})-(\ref{Containment}). We convert this constraint satisfaction problem into an unconstrained optimization problem by introducing a penalty function approach as follows.

\begin{figure*}[!htbp]
\centering
\includegraphics [width=2.0in]{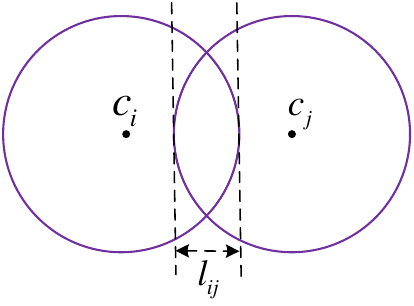}
\caption{Degree of constraint violation ($l_{ij}$) on the distance between two points $c_i$ and $c_j$, where the radius of two circles is $\frac{D}{2}$.}
\label{overlap_betweencircles}
\end{figure*}

Given an allowed minimum distance $D$ between dispersion points and a multiply-connected region $\Omega$ consisting of a polygonal container $C$,  $K$ ($K\ge 0$) polygonal holes $\{U_{1}, U_{2}, \dots, U_{K}\}$, and a candidate solution $X=(x_1,y_1,\dots, x_p,y_p)$ (i.e., Cartesian coordinates of $p$ points in $\Omega$), the quality of $X$ is measured with the following objective function:
\begin{equation}\label{Total_Energy}
E_D(X) = \sum_{i=1}^{p-1}\sum_{j=i+1}^{p} l_{ij}^2 + \alpha \times (\sum_{i=1}^{p} O_{i0} + \sum_{i=1}^{p}\sum_{k=1}^{K} O_{ik})
\end{equation}
with
\begin{equation}\label{EBetweenCircles}
l_{ij} = max\{0, D - \sqrt{(x_i-x_j)^2 + (y_i - y_j)^2}\}
\end{equation}
where $\alpha$ is a penalty factor and $l_{ij}^{2}$ represents the degree of constraint violation of the distance between two dispersion points $c_i$ and $c_j$, and the goal of including $l_{ij}^{2}$ is to force the dispersion points $c_i$ and $c_j$ to be separated by a distance of at least $D$. Equivalently, $l_{ij}$ is the overlap depth between two equal circles with a radius of $\frac{D}{2}$, located respectively at the points $c_i$ and $c_j$ (see Fig. \ref{overlap_betweencircles} for an illustration). Additionally, $O_{i0}$ is the degree of violation of the containment constraint of point $c_i$ and the distance constraint between $c_i$ and the boundary of container $C$, and $O_{ik}$ is the degree of constraint violation of the distance between $c_i$ and the $k$-th hole $U_{k}$. Thus,  $E_D(X)$ measures the total degree of constraint violation of the given solution $X$, and $E_D(X) = 0$ means that $X$ is a feasible solution, i.e., all $p$ points are contained in the container $C$, the minimum distance between the dispersion points is no less than the given $D$ value, and the minimum distance from a dispersion point $c_i$ ($1\le i \le N$) to all boundaries of the container and holes is no less than $D_b$.

$O_{i0}$ and $O_{ik}$ are two key components of our optimization model that are more complicated to understand and formulate than $l_{ij}$ due to the non-convexity of the container $C$ and the hole $U_{k}$. Therefore, we focus on these components in the following two subsections.

\subsubsection{Penalty term $O_{i0}$ between a dispersion point and the polygonal container.}
\label{OverlapContainer}

The penalty term $O_{i0}$ aims to force the dispersion point $c_i$ into the container and keep it at least a distance of $D_b$ from the container  boundary via repulsion and attraction between the dispersion point and this boundary. To formulate the penalty $O_{i0}$, we first give the following two definitions.

\begin{figure*}[!htbp]
\centering
\subfigure[$min_{v\in H_0^{i}} \{ d(c_i,v)\} < min_{v\in V_0} \{ d(c_i,v)\}$] {\includegraphics[width=2.8in]{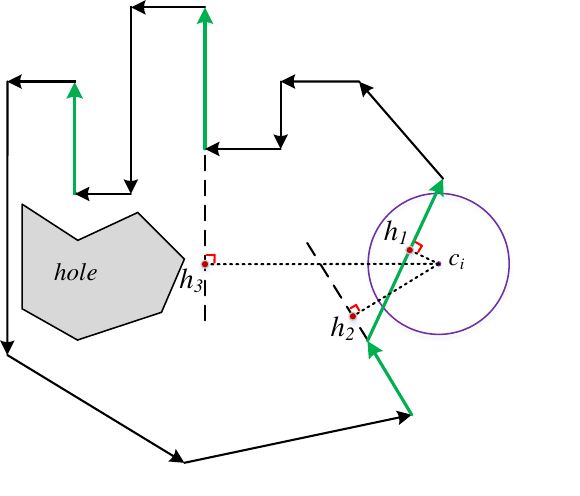}}
\subfigure[$min_{v\in H_0^{i}} \{ d(c_i,v)\} \ge min_{v\in V_0} \{ d(c_i,v)\}$]{\includegraphics[width=2.8in]{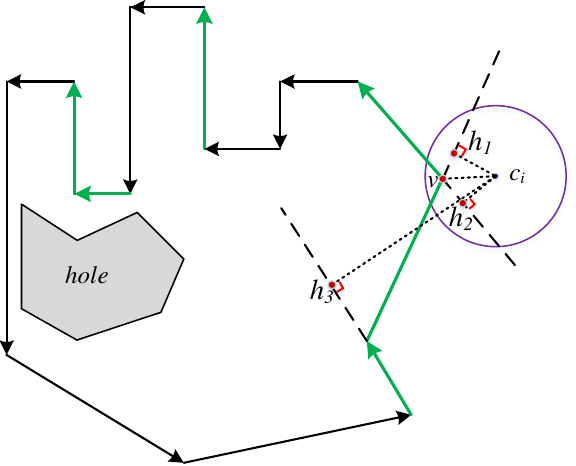}}
\caption{The active edges and some representative foot points for the case that the dispersion point $c_i$ lies outside the container, where the active edges are indicated in green and the circle is $B(c_i,D_b)$.}
\label{ContainerOut}
\end{figure*}

\begin{itemize}
  \item \textbf{Definition 1: Active edges of polygonal container}. Given a polygonal container $C=(V_0,E_0)$ and a dispersion point $c_i$ with coordinates $(x_i,y_i)$, where $V_0$ is the set of vertices of $C$ and $E_0$ is the set of edges on the boundary of $C$, an edge $e_l$ ($\in E_0$) is called an \textit{active edge} with respect to the given point $c_i$ if it satisfies the following condition.

      Moving along the container boundary in a counter-clockwise direction, the current edge $e_l$ ($1 \le l \le |E_0|$) is identified as an active edge in the following two cases: (1) the given point $c_i$ is contained in $C$ and lies on its left-hand side (i.e., using the left-hand rule), (2) the given point $c_i$ is not contained in $C$ and lies on its right-hand side (i.e., using the right-hand rule).

      To illustrate, all active container edges are indicated in green in Figs. \ref{ContainerIn} and \ref{ContainerOut}. Clearly, the answer to the question of whether an edge of the polygonal container is active or not depends on the current position of the dispersion point $c_i$.

  \item \textbf{Definition 2: Active foot point on the container boundary}. Given a dispersion point $c_i$ with coordinates $(x_i,y_i)$ and active edge $e$ on the container boundary, if the foot $h$ of the perpendicular from $c_i$ to the line containing edge $e$ lies on this edge, then the point $h$ is called an \textit{active foot point} with respect to the given dispersion point $c_i$, and a non-active foot point otherwise. The set of all active foot points with respect to the given dispersion point $c_i$ and the container $C$ is denoted by $H_0^i$.

      Figs. \ref{ContainerIn} and \ref{ContainerOut} provide some representative foot points, including active foot points and non-active foot points. Taking the subfigure (a) of Fig. \ref{ContainerIn} as an example, the points $h_1$ and $h_2$ are active foot points, and $h_3$ is a non-active foot point since it lies on the extending line of an edge.
\end{itemize}

Using the set $H_0^{i}$ of active foot points, we formulate the penalty term $O_{i0}$ as follows.
\begin{equation}\label{OverlapWithContainer}
O_{i0} = \begin{cases}
		\gamma \times (D_b + min\{d(c_i,v) : v\in V_0 \cup H_0^{i}\})^2,& c_i \notin C \\
		\sum_{v \in V_0 \cup H_0^{i}} (max\{0,D_b - d(c_i,v)\})^2,& c_i \in C
		\end {cases}
\end{equation}
where $\gamma$ is a penalty factor, $d(c_i,v)$ denotes the Euclidean distance between the point $c_i$ and point $v \in V_0 \cup H_0^{i}$, $V_0$ denotes the set of vertices of the polygonal container, and $D_b$ is the allowed minimum distance from the dispersion point $c_i$ to the container boundaries and can be written as $D_b = \lambda \times  D$ with $\lambda \in [0,0.5]$, where $\lambda$ is a predetermined coefficient used to control the value of $D_b$.

One can observe from Eq. (\ref{OverlapWithContainer}) that the value of $O_{i0}$ is calculated according to whether or not the point $c_i$ is contained in the container $C$. The fundamental principles behind the equation can be explained as follows:
\begin{figure*}[!htbp]
\centering
\subfigure[$H_0^{i} \cap B(c_i,D_b) \neq \emptyset$] {\includegraphics[width=2.5in]{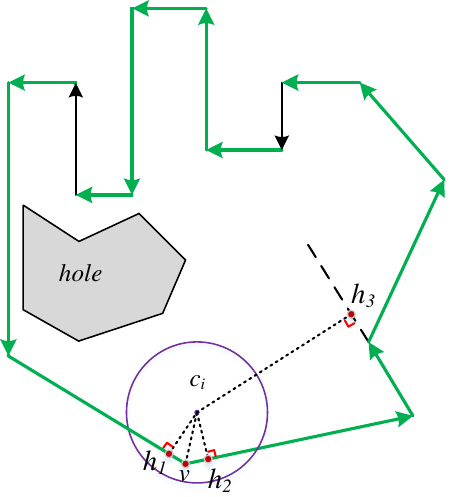}}
\subfigure[$H_0^{i} \cap B(c_i,D_b)   = \emptyset$]{\includegraphics[width=2.5in]{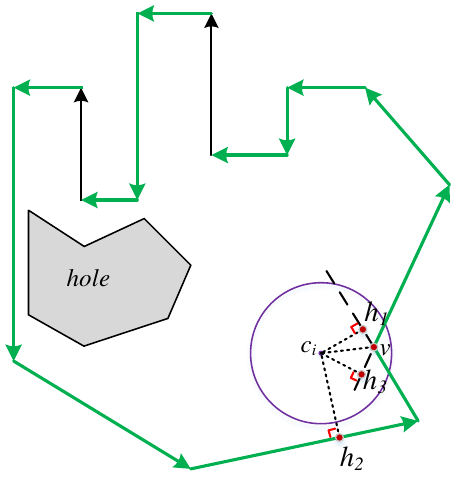}}
\caption{The active edges and some representative foot points for the case that the dispersion point $c_i$ lies in the container, where the active edges are indicated in green and the circle is $B(c_i,D_b)$.}
\label{ContainerIn}
\end{figure*}

First, when the dispersion point $c_i$ with coordinates $(x_i,y_i)$ is contained in $C$ (i.e., $c_i \in C$), we need to identify the degree of constraint violation of the distance between the point $c_i$ and the container boundary to keep the point $c_i$ at least a distance of $D_b$ from the container boundary, which requires consideration of two situations. In the first situation, there exist one or several active foot points in the vicinity $B(c_i,D_b)$ of $c_i$ (i.e., $H_0^{i} \cap B(c_i,D_b) \neq \emptyset$), where $B(c_i,D_b) = \{c: ||c_i -c|| \le D_b\}$. In this case, all points in $(V_0 \cup H_0^{i})\cap B(c_i,D_b)$ are used to generate the penalty term $O_{i0}$, thus producing a repulsive force $(-\frac{\partial O_{i0}}{\partial x_i}, -\frac{\partial O_{i0}}{\partial y_i} )$ to the dispersion point $c_i$, driving it to move a distance of at least $D_b$ away from the boundary. Fig. \ref{ContainerIn}(a) gives a representative example for this situation. On the contrary, if there does not exist any active foot point in the vicinity of $c_i$ (i.e., $H_0^{i} \cap B(c_i,D_b) = \emptyset$), all vertices of container $C$ in the range of $B(c_i,D_b)$ are used to generate $O_{i0}$ and thus produce a repulsive force to $c_i$. Otherwise, it is likely that the dispersion point $c_i$ leaves the container by passing beyond one of its vertices during the optimization process due to the repulsive forces from other dispersion points, which is one of main reasons why the set $V_0$ of container vertices is involved in $O_{i0}$. Fig. \ref{ContainerIn}(b) gives an example of this situation. Formally, $O_{i0}$ can be written as $O_{i0}= \sum_{v \in V_0 \cup H_0^{i}} (max\{0,D_b - d(c_i,v)\})^2$ for these two situations.

Second, when the point $c_i$ lies outside the container $C$ (i.e., $c_i \notin C$), two situations need to be considered to handle the containment constraint of point $c_i$, to force the point $c_i$ to enter the container $C$. In the first situation, the minimum distance between $c_i$ and the set $H_0^{i}$ is smaller than the minimum distance between $c_i$ and the set $V_0$ (i.e., $min_{v\in H_0^{i}} \{ d(c_i,v)\} < min_{v\in V_0} \{ d(c_i,v)\}$). Then the foot point $v \in H_0^{i}$ closest to $c_i$ is used to generate the penalty term $O_{i0}$ and thus produce an attractive force to induce $c_i$ to gradually approach the boundary of container. Fig. \ref{ContainerOut}(a) provides an example for this situation. In the second situation, $min_{v\in H_0^{i}} \{ d(c_i,v)\} \ge min_{v\in V_0} \{ d(c_i,v)\}$, and then the vertex $v$ of $V_0$ closest to $c_i$ is used to generate the penalty term $O_{i0}$ and thus produces an attractive force $(-\frac{\partial O_{i0}}{\partial x_i}, -\frac{\partial O_{i0}}{\partial y_i} )$ between the points $v$ and $c_i$. Without these adjustments, the dispersion point $c_i$ will always lie outside of container during the optimization due to the lack of an attractive force, which is another reason why the set $V_0$ of vertices is involved in $O_{i0}$. Fig. \ref{ContainerOut}(b) provides an example for this situation. In summary, when $c_i \notin C$ holds, the point in $V_0 \cup H_0^{i}$ closest to $c_i$ is used to generate $O_{i0}$. Formally, the penalty term $O_{i0}$ can be uniformly written as $O_{i0}= \gamma \times (D_b + min\{d(c_i,v): v\in V_0 \cup H_0^{i}\})^2$ for these two situations.

Thus, the inclusion of $O_{i0}$ in $E_D(X)$ ensures that the dispersion point $c_i$ is contained in the container and the minimum distance between $c_i$ and boundary of container is no less than $D_b$ for any candidate solution $X$ with $E_D(X)=0$.

\subsubsection{Penalty term $O_{ik}$ between a dispersion point and a polygonal hole.}
\label{OverlapHole}

The penalty term $O_{ik}$ aims to force the dispersion point $c_i$ to leave from the hole $U_k$ and keep at least a distance $D_b$ from the hole boundary. To formulate $O_{ik}$, we first give two definitions.

\begin{figure*}[!htbp]
\centering
\subfigure[$min_{v\in V_k}\{d(c_i,v)\} < min_{v \in H_k^{i}} \{d(c_i,v)\}$] {\includegraphics[width=2.7in]{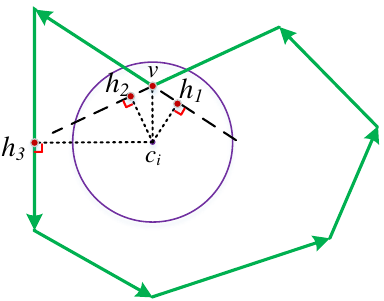}}
\subfigure[$min_{v\in V_k}\{d(c_i,v)\} \ge min_{v \in H_k^{i}} \{d(c_i,v)\}$]{\includegraphics[width=2.7in]{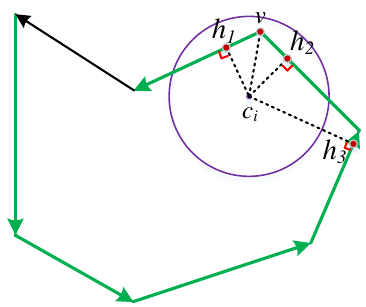}}
\caption{The active edges and some representative foot points for the case that the dispersion point $c_i$ lies in the hole $U_{k}$, where the active edges are indicated in green and the circle is $B(c_i,D_b)$.}
\label{HoleIn}

\end{figure*}

\begin{figure*}[!htbp]
\centering
\subfigure[$H_k^{i} \cap B(c_i,D_b)  \neq \emptyset$ ]{\includegraphics[width=2.5in]{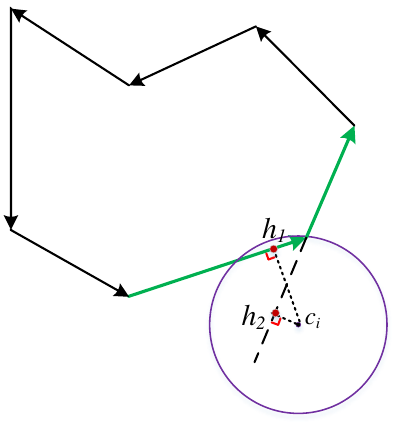}}
\subfigure[$H_k^{i} \cap B(c_i,D_b)  =  \emptyset$]{\includegraphics[width=2.5in]{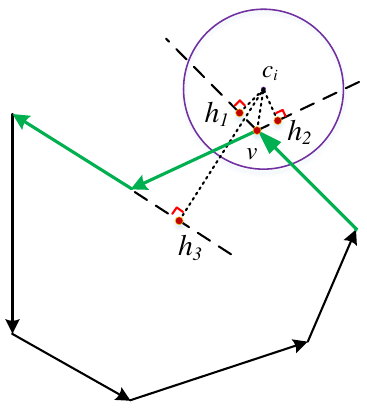}}
\caption{The active edges and some representative foot points for the case that the dispersion point $c_i$ lies outside the hole $hole_{k}$, where the active edges are indicated in green and the circle is $B(c_i,D_b)$.}
\label{HoleOut}
\end{figure*}

\begin{itemize}
  \item \textbf{Definition 3: Active edges of a polygonal hole}. Given a polygonal hole $U_k=(V_k,E_k)$ ($k\ge 1$) and a dispersion point $c_i$ with coordinates $(x_i,y_i)$, an edge $e$ ($\in E_k$) is called an \textit{active edge} of $U_k$ with respect to $c_i$ according to two situations. Specifically, when moving along the boundary of hole $U_k$ in a counter-clockwise direction, the edge $e$ is identified as an active edge of $U_k$ if (1) the given point $c_i$ is contained in $U_k$ and lies on its right-hand side (i.e., using the right-hand rule), or (2) the given point $c_i$ is not contained in $U_k$ and lies on its left-hand side (i.e., using the left-hand rule).

      For example, all active edges of a polygonal hole are indicated in green in Figs. \ref{HoleIn} and \ref{HoleOut}. Clearly, the answer to the question of whether an edge of a polygonal hole is active or not depends on the current position of the dispersion point $c_i$.

  \item \textbf{Definition 4: Active foot point on the boundary of a polygonal hole}. Given a dispersion point $c_i$ with coordinates $(x_i,y_i)$ and an active edge $e$ of the hole $U_k$, if the foot $h$ of perpendicular from $c_i$ to the line containing the edge $e$ lies on $e$, then $h$ is called an \textit{active foot point} of $U_k$ with respect to the dispersion point $c_i$. The set of all active foot points of $U_k$ with respect to the given dispersion point $c_i$ is denoted by $H_k^i$.

      Figs. \ref{HoleIn} and \ref{HoleOut} give some representative foot points, including active foot points and non-active foot points. In Fig. \ref{HoleIn}(a), the point $h_3$ is an active foot point, and the points $h_1$ and $h_2$ are non-active foot points since they lie on the extended line of the corresponding edge.
\end{itemize}

Using the set $H_k^{i}$ of active foot points, the penalty term $O_{ki}$ can be formulated as follows:
\begin{equation}\label{OverlapWithHoles}
O_{ik} = \begin{cases}
		\gamma \times (D_b + min\{d(c_i,v) : v\in V_k \cup H_k^{i}\})^2,& c_i \in U_k \\
		\sum_{v \in V_k \cup H_k^{i}} (max\{0,D_b - d(c_i,v)\})^2,& c_i \notin U_k
		\end {cases}
\end{equation}
where $\gamma$ is a penalty factor, $d(c_i,v)$ represents the distance between point $c_i$ and $v$, $V_k$ is the set of vertices of $U_k$, $H_k^{i}$ is the set of active foot points of $U_k$ with respect to the given dispersion point $c_i$, and $D_b$ represents the allowed minimum distance from $c_i$ to the boundaries of holes and container.


Equation (\ref{OverlapWithHoles}) indicates that the penalty term $O_{ik}$ is calculated depending on whether the dispersion point $c_i$ lies in the hole $U_k$ or not. First, if $c_i \in U_k$, the point in $V_k \cup H_k^{i}$ closest to $c_i$ is used to generate $O_{ik}$ and the dispersion point $c_i$ will receive an attractive force $(-\frac{\partial O_{ik}}{\partial x_i}, -\frac{\partial O_{ik}}{\partial y_i} )$ from this point, thus guiding it to leave the hole $U_k$. Thus, $O_{ik}$ can be written as $\gamma \times (D_b + min\{d(c_i,v) : v\in V_k \cup H_k^{i}\})^2$. Fig. \ref{HoleIn} gives two examples for the following two situations: (1) the closest point belongs to $V_k$ (subfigure (a)), (2) the closest point belongs to $H_k^{i}$ (subfigure (b)).

Second, if the dispersion point $c_i$ lies on the outside of hole $U_k$ (i.e., $c_i \notin U_k$), all the points in $(V_k \cup H_k^{i})\cap B(c_i,D_b)$ will provide a positive contribution to the penalty term $O_{ik}$ and thus produce a repulsive resultant force $(-\frac{\partial O_{ik}}{\partial x_i}, -\frac{\partial O_{ik}}{\partial y_i} )$ on the dispersion point $c_i$, making it away from the hole boundary by a distance of at least $D_b$. If there is no active foot point of $U_k$ in the vicinity $B(c_i, D_b)$ of $c_i$, then the vertices of the hole $U_k$ in $B(c_i, D_b)$ are used to generate the penalty term $O_{ik}$. An example of this situation can be found in subfigure (b) of Fig. \ref{HoleOut}. Thus, $O_{ik}$ can be written as $O_{ik}= \sum_{v \in V_k \cup H_k^{i}} (max\{0,D_b - d(c_i,v)\})^2$.

It is worth noting that, according to Eqs. (\ref{OverlapWithContainer}) and (\ref{OverlapWithHoles}), the point‐in‐polygon test \citep{wang20052d} is a basic subroutine that will be frequently called during the optimization process for any algorithm employing this optimization model.

\subsection{The parameters and discussion of the proposed optimization model}
\label{Model_parameters}

The proposed optimization model involves three parameters  $\alpha$, $\gamma$, and $\lambda$, where $\alpha$ is used in Eq. (\ref{Total_Energy}), and $\gamma$ is used in Eqs. (\ref{OverlapWithContainer}) and (\ref{OverlapWithHoles}) to quantify the penalty degree of constraint violation. Parameter $\lambda$ is used to control the value of $D_b$ as $D_b = \lambda \times D$, where $D$ is the allowed minimum distance between dispersion points and $D_b$ is the allowed minimum distance between a dispersion point and the boundaries of the container and holes, where a larger $\lambda$ value means a larger distance constraint, and vice versa.

The optimization model in Eq. (\ref{Total_Energy}) is uniformly used to formulate the continuous $p$-dispersion problems with and without a boundary constraint. To distinguish these two problems, we use respectively different parameter settings as follows.
\begin{itemize}
   \item For the continuous $p$-dispersion problem with a boundary constraint (i.e., the equal circle packing problem), the parameters $\alpha$, $\gamma$, and $\lambda$ are respectively set to  1.0, 2.0 and 0.5, and the resulting model is equivalent to the equal circle packing problem.
   \item For the continuous $p$-dispersion problem without a boundary constraint (i.e., the point arrangement problem), the parameters $\alpha$, $\gamma$, and $\lambda$ are respectively set to 3000, 1.0 and $10^{-3}$. We set $\lambda$ to a very small number (i.e., $10^{-3}$), instead of 0, as our preliminary experiments disclosed that $\lambda=0$ results in an ill-conditioned objective function for the problems with a non-convex region, significantly increasing the difficulty of local optimization.
\end{itemize}

One notices that the proposed optimization model is differentiable except for the points on the boundaries of the container and holes, making it possible to apply popular continuous solvers (e.g., the quasi-Newton method) as a local optimization method to reach high-precision solutions of the studied problem. Until now, such a dedicated optimization model has been missing for the continuous $p$-dispersion problems with a non-convex multiply-connected region, which explains why existing algorithms in the literature fail to obtain high-precision solutions for these cases. The proposed unified optimization model fills this gap and can be considered as one of the main innovations of the present work.

\section{Tabu Search based Global Optimization for the continuous $p$-dispersion problems}
\label{Method}

Our TSGO algorithm is designed to solve the optimization model described in Section \ref{Models} for the continuous $p$-dispersion problems with and without boundary constraints. The basic idea of the algorithm is to dynamically convert a constrained optimization problem into a series of unconstrained optimization subproblems and then solve them by a tabu search method operating in the continuous solution space.

\subsection{Main framework of the TSGO algorithm}
\label{general_procedure}

\renewcommand{\baselinestretch}{1.2}\huge\normalsize
\begin{algorithm}[h]
\small
\DontPrintSemicolon 
\KwIn{Number of points to be packed ($p$), container $C$ and holes $U_{k}$ ($k=1,2,\dots, K$), maximum time limit ($t_{max}$), packing density of the initial solution ($\rho$) }
\KwOut{The best configuration found ($X^{*}$) and the minimum distance between points ($D^{*}$)}
$Area \gets CalculateArea(C, U_{1}, U_{2}, \dots, U_{K})$\tcc*[r]{Calculate the area of bounded region $\Omega$ as $Area = Area(C) -\sum_{k=1}^{K} Area(U_{k})$}
$D \gets 2\times \sqrt{\frac{\rho \times Area}{N \times \pi}}$ \tcc*[r]{Calculate the initial value of $D$ according to input density $\rho$}
$X \gets RandomSolution(p)$\tcc*[r]{Generate randomly an initial solution $X$}
$X \gets TabuSearch(E_D(X),X,D)$\tcc*[r]{Minimize the function $E_D(X)$ by Algorithm \ref{algo:TS}}
$(X,D)\gets AdjustDistance(X,D)$\tcc*[r]{Adjust configuration $(X,D)$ by Algorithm \ref{algo:adust}}
$D^{*} \gets D$, $X^{*} \gets X$\;
\While{\textit{time()} $\le$ $t_{max}$}
{
  $D \gets D^{*}$\tcc*[r]{Set $D$ to the best value found so far}
  $X \gets RandomSolution(p)$\;
  $(X,D) \gets TabuSearch(E_D(X),X,D)$\;
  \tcc*[l]{$E_D(X) < 10^{-25}$ means that $X$ is a feasible solution for the current $D$}
  \If{$E_D(X) < 10^{-25}$}
  {
    $(X,D)\gets AdjustDistance(X,D)$\;
    \If{$D > D^{*}$}
    {
       $D^{*} \gets D$\;
       $X^{*} \gets X$\tcc*[r]{Save the best solution found}
    }
  }

}
\Return{$(X^{*},D^{*})$}\;
\caption{Main framework of the proposed TSGO algorithm}
\label{algo:TSGO}
\end{algorithm}
\renewcommand{\baselinestretch}{1.0}\huge\normalsize

The TSGO algorithm consists of three main components: a solution initialization procedure, a tabu search method to improve the solution, and a distance adjustment method to maximize the minimum distance between dispersion points while maintaining the feasibility of the solution. The pseudo-code of the algorithm is described in Algorithm \ref{algo:TSGO}, where $X$ represents the current solution, $X^*$ and $D^*$ respectively denote the best solution found so far and the corresponding minimum distance between the dispersion points.

The algorithm works as follows. First, it calculates the area of the bounded region $\Omega$ (line 1), which is the difference between the area of the container and the total area of the $K$ holes. Then, assuming that each dispersion point is a circle with a radius of $\frac{D}{2}$ and that the preestimated packing density of these $p$ circles in the region  $\Omega$ is an input value $\rho$, the algorithm calculates $D = 2\times \sqrt{\frac{\rho \times Area}{N \times \pi}}$, which is the initial estimation of the minimum distance between the $p$ points (line 2).

After that, the algorithm generates an initial solution by distributing uniformly and randomly the $p$ dispersion points in a minimum rectangle containing the region $\Omega$ and then employs the tabu search method to improve the quality of solution by minimizing the function $E_D(X)$ (lines 3-4). The improved solution from tabu search is slightly adjusted to maximize the minimum distance $D$ between $p$ points,  maintaining the feasibility of the solution (line 5), and then the resulting solution is saved as $X^{*}$.

Then, the algorithm enters a `while' loop to iterate until the time limit ($t_{max}$) is reached and finally returns the best solution $X^{*}$ found (lines 7-19). At each iteration, based on the best value of $D$ found so far, an initial solution is randomly generated and the tabu search method is used to improve its quality by minimizing $E_D(X)$ (lines 8-10). The distance adjustment procedure ($AdjustDistance(\cdot)$) is used to maximize the minimum distance $D$ between points once a feasible solution is obtained by the tabu search method (line 13). Moreover, the value of $D^{*}$ is updated each time an improving value is found (lines 14-17).

\subsection{Tabu search method}
\label{subsec_TS}

Tabu search (TS) \citep{glover1998tabu} is a popular metaheuristic for combinatorial optimization that has also been adapted to solve continuous global optimization problems \citep{chelouah2000tabu}. The present TS method exploits the solution space utilizing a neighborhood composed of nearby local minima. The general procedure, neighborhood structure and tabu list management strategy of our tabu search method are described in the following subsections.

\subsubsection{Tabu search for the continuous $p$-dispersion problems.}
\label{subsec_tabu_search}

The TS method of our TSGO algorithm is depicted in Algorithm \ref{algo:TS}, where $X$ and $X^{b}$ respectively represent the current solution and the best solution found so far, and $X_{neighbor}^{best}$ represents the best individual in the current neighborhood $N_{ins}(X)$ in terms of the objective value. Starting from an input solution $X_0$ and the initialization of the tabu list $T$ (lines 2-5), the TS method performs successive iterations until a feasible solution is found or the best solution $X^{b}$ cannot be improved during the last $\beta_{max}$ consecutive iterations (lines 6-31), where $\beta_{max}$ is a parameter called the search depth. At each iteration, the candidate set $P$ of high-energy dispersion points and the candidate set $V$ of low-energy vacancy sites are first constructed for the current solution $X$ (lines 7-8). Then, the vacancy-based insertion neighborhood $N_{ins}(X)$ of the current solution with respect to $P$ and $V$ is evaluated (lines 10-21) and a best non-tabu solution in $N_{ins}(X)$ is selected to replace the current solution $X$ (line 22). After that, the tabu list $T$ is updated and the current solution $X$ is further improved by a very short run of the monotonic basin-hopping (MBH) method, which is proposed in \citep{leary2000global} for an intensified search (lines 23-24).

The MBH method is a simple iterated improvement procedure, which repeats a perturbation operator followed by a local optimization until the current solution cannot be further improved for $\theta_{max}$ iterations, where $\theta_{max}$ is a parameter called its search depth. At each iteration of the MBH method, the current solution is first perturbed by shifting each solution coordinate of the solution in a given interval $[-\Delta, \Delta]$, where $\Delta = 0.4 \times D$, and then is locally improved by the L-BFGS method \citep{liu1989limited} equipped with an efficient line search method \citep{hager2005new}. The new solution is accepted as the current solution if and only if it is superior to the current solution.

\renewcommand{\baselinestretch}{1.2}\huge\normalsize
\small
\begin{algorithm}[h]
\DontPrintSemicolon 
\textbf{Function} \emph{TabuSearch}\;
\KwIn{Input solution $X_0$, depth of tabu search ($\beta_{max}$), parameter $Q$.}

\KwOut{The best solution found ($X^{b}$)}
$X \gets Local\_Optimization(X_0)$\tcc*[r]{Minimize locally function $E_D(\cdot)$ starting from $X$}
$X^{b} \gets X$\tcc*[r]{$X^{b}$ denotes the best solution found by the current tabu search}
$NoImprove \gets 0$\;
Initialize tabu list $T$\;
\While{($NoImprove$ $\le$ $\beta_{max}$) $\land$ ($E_{D}(X^{b}) > 10^{-25}$) }
{
  Select $Q$ dispersion points (i.e., $P[1:Q]$) with the highest energies from the current solution $X$\;
  Find $Q$ vacancy sites (i.e., $V[1:Q]$)) with the lowest energies in the current solution $X$\;
  \tcc*[l]{Evaluate the insertion neighborhood $N_{ins}(X)$}
  $i_1 \gets rand(Q)$, $j_1 \gets rand(Q)$\tcc*[r]{$rand(Q)$ denotes a random integer in $[1,Q]$}
  $X_{neighbor}^{best} \gets  X \oplus Move(P[i_1],V[j_1])$\tcc*[r]{$X_{neighbor}^{best}$ denotes the best solution in $N_{ins}(X)$}
  \For{$i \gets 1$ \textbf{to} $Q$ }
  {
     \For{$j \gets 1$ \textbf{to} $Q$ }
     {
        $X_{neighbor} \gets X \oplus Move(P[i],V[j])$\;
        $X_{neighbor} \gets Local\_Optimization(X_{neighbor})$\tcc*[r]{Minimize $E_D(X)$ by L-BFGS}
        \If{$Move(P[i],V[j])$ is not forbidden $\land$ $E_D(X_{neighbor})< E_D(X_{neighbor}^{best})$}
        {
           $X_{neighbor}^{best} \gets X_{neighbor}$\;
           $I \gets P[i]$\;
        }
     }
  }

  $X \gets X_{neighbor}^{best}$\tcc*[r]{Update the current solution}
  Update tabu list $T$ by $I$\;
  $X \gets MBH(X)$ \tcc*[r]{Improve $X$ by a short run of MBH}
  \uIf{$E_{D}(X) < E_D(X^{b})$}
  {
      $X^{b} \gets X$ \tcc*[r]{Save the best solution found}
      $NoImprove \gets 0$\;
  }
  \Else
  {
      $NoImprove \gets NoImprove + 1$\;
  }
}
\Return{$X^{b}$}\;
\caption{Tabu Search method for the continuous $p$-dispersion problems}
\label{algo:TS}
\end{algorithm}
\renewcommand{\baselinestretch}{1.0}\huge\normalsize

\subsubsection{Neighborhood structure.}
\label{subsec_neighborhood}

\begin{figure*}[!htbp]
\centering
\includegraphics [width=2.2in]{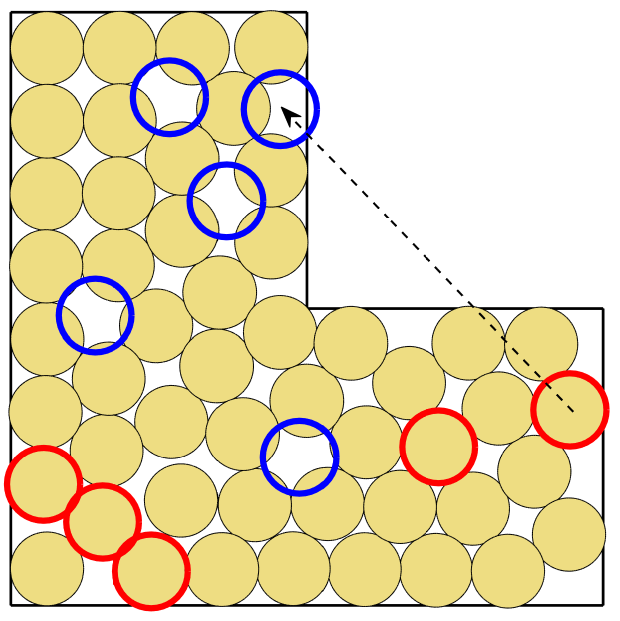}
\caption{The neighborhood move used in the tabu search method.}
\label{N1}
\end{figure*}

The neighborhood structure is a highly important component of tabu search algorithms and is usually applied in the setting of discrete optimization. However, for the continuous $p$-dispersion problems, an efficient neighborhood structure is not easy to design. In this study, we devise a vacancy-based insertion neighborhood $N_{ins}(X)$ based on the potential energies of the $p$ dispersion points and the vacancy sites in the current configuration $X$, which is defined by the insertion operator that moves a high-energy dispersion point to a low-energy vacancy site.

Taking the continuous $p$-dispersion problem with a boundary constraint (i.e., the equal circle packing problem) as an example, we describe the vacancy-based insertion neighborhood $N_{ins}(X)$ as follows. First, the potential energy $E(c_i)$ of each circle $c_i$ in $X$, which measures the degree of constraint violation of $c_i$ in $X$, is calculated by using Eq. (\ref{EnergySinglePoint}), and the first $Q$ highest-energy circles constitute the candidate set $P$, where $Q$ is a parameter.
\begin{equation}\label{EnergySinglePoint}
E(c_i) = \sum_{j=1,j\neq i}^{p} l_{ij}^2 + \alpha \times (O_{i0} + \sum_{k=1}^{K} O_{ik})
\end{equation}
and $l_{ij}$, $O_{i0}$ and $O_{ik}$ are defined in Eqs. (\ref{EBetweenCircles})--(\ref{OverlapWithHoles}).


Then, to construct the candidate set $V$ of low-energy vacancy sits, the algorithm performs $5p$ random detections. Specifically, for each random detection, an additional circle called the detecting circle is first placed randomly in the region $\Omega$ and then its position is locally optimized via minimizing the potential function Eq. (\ref{EnergySingleVacancy}), which measures the degree of constraint violation when a circle is placed at its current position, by means of the L-BFGS method \citep{liu1989limited}. After that, the local minimum solution returned by L-BFGS method is identified as the vacancy site detected. Finally, the first $Q$ lowest-energy vacancy sites found are used to generate the candidate set $V$. It is worth noting that the minimization procedure for the potential energy function (\ref{EnergySingleVacancy}) is very fast and the average computational time is about $10^{-4}$ seconds on our computer when $X$ contains 100 circles (or dispersion points).
\begin{equation}\label{EnergySingleVacancy}
E(v_i) = \sum_{j=1}^{p} l_{ij}^2 + \alpha \times (O_{i0} + \sum_{k=1}^{K} O_{ik})
\end{equation}

With the help of the candidate sets $P$ and $V$, the vacancy-based insertion neighborhood $N_{ins}(X)$ of the current solution $X$ can be easily defined as follows:
\begin{equation}\label{Neighborhood}
N_{ins}(X) = \{LS(X \oplus Move(P[i],V[j])): 1 \le i,j \le Q \}\;
\end{equation}
where $Move(P[i],V[j])$ denotes the insertion operation, which moves the circle (or dispersion point) $P[i]$ from its current position to the vacancy site $V[j]$, and $LS(\cdot)$ denotes a local optimization procedure such as L-BFGS. Clearly, the size of the neighborhood $N_{ins}(X)$ equals $Q^{2}$. Figure \ref{N1} illustrates an insertion operation, where the first $Q$ ($Q=5$) highest-energy circles are indicated in red and the first $Q$ lowest-energy vacancy sites obtained by the detection procedure are indicated with a blue circle.

\subsubsection{Tabu list management.}
\label{subsec_tabu_list}

When a neighborhood move $Move(c_i,v_j)$ is performed, the dispersion point $c_i$ is forbidden to move by the insertion operation for the next $tt(Iter)$ iterations, where $Iter$ is the  iteration number. The tabu tenure $tt(Iter)$ is simply determined as $tt(Iter) = 5 + rand(0,5)$, where $rand(0,5)$ is a random number in [0,5.

\subsection{Adjustment method for the minimum distance between dispersion points}
\label{subsec_distance_adjust}

\renewcommand{\baselinestretch}{1.2}\huge\normalsize
\begin{algorithm}[htbp]
\small
\DontPrintSemicolon 
\textbf{Function} \emph{AdjustDistance}\;
\KwIn{Input solution $(X_0,D_0)$, maximum number of iterations $K$ ($=15$)}
\KwOut{Feasible local optimum solution $(X,D)$}
$X \gets X_0$, $D \gets D_0$, $\mu \gets 10$\;
\For{$i \gets 1$ \textbf{to} $K$ }
{
   $(X,D)$  $\gets$  $Local\_Optimization(\Phi_{\mu},X,D)$ \tcc*[r]{Minimize $\Phi_{\mu}(X,D)$ by L-BFGS method}
   $\mu \gets 5 \times \mu$\;
}
\Return{$(X,D)$}\;
\caption{Adjustment method for the minimum distance $D$ between points}
\label{algo:adust}
\end{algorithm}
\renewcommand{\baselinestretch}{1.0}\huge\normalsize

Given an input solution $X$ consisting of $p$ dispersion points in region $\Omega$ and an estimated value $D_0$ for the minimum distance between the dispersion points, the distance adjustment procedure adjusts slightly the coordinates of the points of $X$, to maximize the minimum distance $D$ between the dispersion points while preserving the feasibility of the resulting solution. Theoretically, this is equivalent to finding a locally optimal solution closest to the input solution $X$ for the constrained optimization problem defined by Eqs. (\ref{DMax})--(\ref{Containment}). To find a local optimal solution for a constrained optimization problem, we employ the well-known sequential unconstrained minimization technique (SUMT) \citep{fiacco1964computational}. The basic idea of SUMT (Algorithm \ref{algo:adust}) is to progressively convert a constrained optimization problem into a series of unconstrained subproblems, and solves them sequentially until reaching a feasible local minimum solution for the original problem.

The constrained optimization problem corresponding to the continuous $p$-dispersion problem is defined by Eqs. (\ref{DMax})--(\ref{Containment}) and can be progressively converted into the following unconstrained subproblems.
\begin{equation}\label{Penalty}
Minimize \quad \Phi_\mu (X,D) = -D^2 + \mu \times E(X,D)
\end{equation}
where $X$ denotes a candidate solution, $\mu$ is a penalty factor whose each value defines an unconstrained optimization function $\Phi_{\mu}(X,D)$, $D$ is a variable representing the allowed minimum distance between $p$ dispersion points, and the term $E(X,D)$ measures the degree of constraint violation in $X$ and can be formulated as follows:
\begin{equation}\label{Energy}
E(X,D) = \sum_{i=1}^{p-1}\sum_{j=i+1}^{p} l_{ij}^2 + \alpha \times (\sum_{i=1}^{p} O_{i0} + \sum_{i=1}^{p}\sum_{k=1}^{K} O_{ik})
\end{equation}
where $l_{ij}$, $O_{i0}$ and $O_{ik}$ are respectively defined in Eqs. (\ref{EBetweenCircles})--(\ref{OverlapWithHoles}). It is worth noting that $E(X,D)$ defined by Eq. (\ref{Energy}) differs from $E_D(X)$ defined by Eq. (\ref{Total_Energy}) in Section \ref{Model1}, since $E(X,D)$ and $E_D(X)$ respectively contain $2N+1$ and $2N$ variables, and $D$ denotes a constant and a variable respectively in the functions $E_D(X)$ and $E(X,D)$.


\section{Computational Experiments and Assessments}
\label{Experimental_results}

In this section, we evaluate the performance of the TSGO algorithm by performing extensive experiments on benchmark instances and comparing it with state-of-the-art algorithms from the literature. 

\subsection{Benchmark Instances, Parameter Settings and Experimental Protocol}
\label{exprimental_Protocol}

\renewcommand{\baselinestretch}{1.2}\large\normalsize
\begin{table}[!htbp]\centering
\begin{scriptsize}
\caption{Settings of important parameters}
\label{Parameter_Settings}
\begin{tabular}{p{1.2cm}p{1.2cm}p{5.5cm}p{0.01cm}p{2.0cm}p{0.01cm}}
\hline
Parameters           &Section                      & Description                                  & &  Values           \\
\hline
$\rho$               & \ref{general_procedure}     & packing density of initial solution          & &  $\{0.85,1.4\}$   \\
$Q$                  & \ref{subsec_neighborhood}   & parameter for constructing neighborhood      & &  3                \\
$\beta_{max}$        & \ref{subsec_tabu_search}    & search depth of tabu search                  & &  $\{5, 50\}$      \\
$\theta_{max}$       & \ref{subsec_tabu_search}    & search depth of MBH                          & &  $\{10, 15\}$     \\
\hline
\end{tabular}
\end{scriptsize}
\end{table}
\renewcommand{\baselinestretch}{1.0}\large\normalsize

For the equal circle packing problem, the test bed used includes 12 small-sized instances taken from \citep{dai2021repulsion} and Erich's packing center (\url{https://erich-friedman.github.io/packing/cirinl/}), together with 105 instances that we generated by constructing various regions for various problem sizes, where two regions are taken from the website \url{https://people.sc.fsu.edu/~jburkardt/datasets/polygon/polygon.html}. For the point arrangement problem, the test bed includes 75 instances that we generated by setting various problem sizes for these constructed regions.


The default settings and descriptions of TSGO's parameters are given in Table \ref{Parameter_Settings}, where the values were obtained from preliminary experiments. The parameter $\rho$ was set to 0.85 and 1.4, the parameter $\theta_{max}$ was set to 10 and 15, and the parameter $\beta_{max}$ was set to 50 and 5 for the equal circle packing and point arrangement problems, respectively. The algorithm was implemented in C++ and was run on a computer with an Intel(R) Xeon (R) Platinum 9242 CPU (2.3 GHz). To evaluate the average performance, the TSGO algorithm as well as the main reference algorithms were run 10 times with different random seeds for each tested instance. The stopping criterion of the TSGO algorithm and its main reference algorithms is the maximum time limit $t_{max}$ that varies according to the instance size $p$. For the CpDP problem with boundary constraints, $t_{max}$ was set to 50 seconds for the small instances with $p\le50$, and set to 1, 2 and 4 hours respectively for the large instances with $p=100$, 150 and 200. For the CpDP problem without boundary constraint, $t_{max}$ was to 1, 4 and 8 hours respectively for the instances with $p=50$, 100 and 150, since the local optimization method has a slow convergence for the potential functions employed.

\subsection{Computational results and comparison on the equal circle packing problem}
\label{resutls_withBC}

This section assesses the performance of the TSGO algorithm on the CpDP problem with boundary constraints, which corresponds to the equal circle packing problem. For the continuous $p$-dispersion problems in a simply-connected non-convex region, the repulsion-based dispersion algorithm (RBDA) recently proposed in \citep{dai2021repulsion} can be regarded as the state-of-the-art algorithm and consequently the objective of our first experiment is to compare the TSGO algorithm and the RBDA algorithm on popular benchmark instances in the literature, including 10 small-scale instances from Erich's packing center (\url{https://erich-friedman.github.io/packing/cirinl/}) and 2 larger instances from a previous study \citep{dai2021repulsion}. The comparative results are summarized in Table \ref{results_small}. Columns 1-3 show instance features, including the container, the number ($|E|$) of edges of the container, and the number of dispersion points. Column 4 gives the best results of the RBDA algorithm, which were extracted from the literature. The remaining columns of the table give the detailed results of our TSGO algorithm, including the best objective value $R_{best}$ (i.e., the circle radius, which equals $D_{best}/2$) over 10 independent runs, along with the average and the worst objective values $R_{avg}$ and $R_{worst}$, the success rate of hitting the best solution (SR), the standard deviation of the objective values obtained ($\sigma$), and the average computational times in seconds to hit the best solution. The last row ``\#Best" of the table indicates the number of instances for which the corresponding result is the best in terms of the corresponding performance indicator.

\begin{figure*}[h]
\centering
\subfigure[E11H0 for $p=70$]{\includegraphics[width=2.2in]{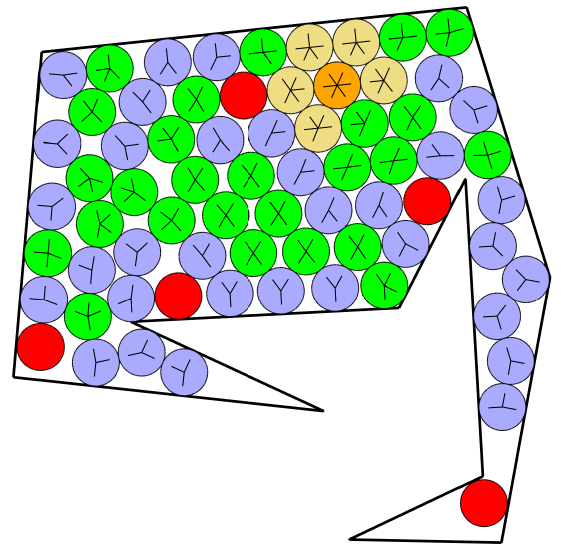}}
\subfigure[E12H0b for $p=100$]{\includegraphics[width=2.2in]{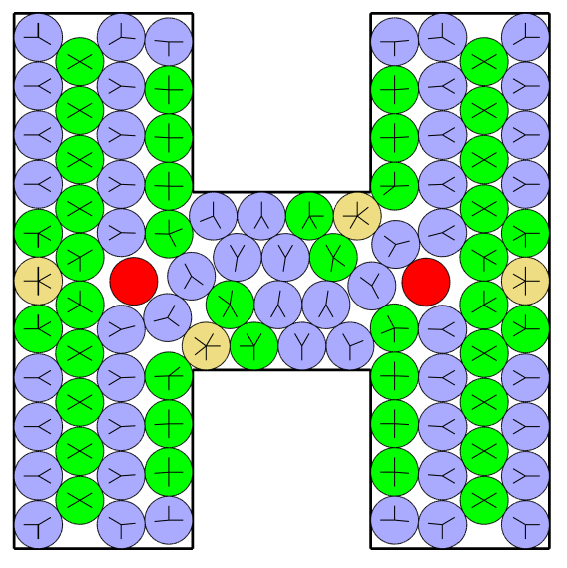}}
\caption{Improved best solutions for two representative instances from the literature. The circles are colored according to their number of neighbors (e.g., all circles with 0 neighbor are colored in red.) }
\label{packing_improved}
\end{figure*}

Table \ref{results_small} shows that our TSGO algorithm significantly outperforms the RBDA algorithm, obtaining a better $R_{best}$ value for all instances tested. Moreover, the success rates of the TSGO algorithm is 100\% for all instances and the computational time is short, especially for the small instances with $p\le 16$.
For very small instances with $p\le 16$, the solutions obtained by both algorithms should have roughly the same geometry, and the differences in the results are caused by the precision of their results. Unlike the RBDA algorithm, TSGO attains a very high precision of $10^{-10}$ for the solutions. This indicates that an appropriate optimization model is very important for the effectiveness of global optimization algorithms for this continuous dispersion problem. For the two larger instances with $p=70$ and 100, the TSGO results are appreciably superior to the RBDA results, and the improved configurations are given in Fig. \ref{packing_improved} to show their differences relative to the previous best-known configurations published in \citep{dai2021repulsion}.

\renewcommand{\baselinestretch}{1.1}\huge\normalsize
\begin{table}[!htbp]\centering
\caption{Comparison of the proposed TSGO algorithm with the best performing algorithm in the literature on 12 popular instances from the literature.} \label{results_small}
\begin{tiny}
\begin{tabular}{p{0.8cm}p{0.4cm}p{0.4cm}p{0.01cm}p{1.0cm}p{0.01cm}p{1.5cm}p{1.5cm}p{1.5cm}p{0.8cm}p{0.5cm}p{0.5cm}p{0.01cm}}
\hline
\multicolumn{3}{c}{\textbf{Instance}} &&  \multicolumn{1}{c}{\textbf{RBDA}} &&  \multicolumn{6}{c}{\textbf{TSGO (This work)}}  &\\
\cline{1-3} \cline{5-5} \cline{7-12}
\textbf{Container} & $|E|$ & $p$  && \centering{\textbf{$R_{best}$} }&&  \centering{\textbf{$R_{best}$}} &  \centering{\textbf{$R_{avg}$}} & \centering{\textbf{$R_{worst}$}}& \centering{\textbf{$SR$}} &\centering{\textbf{$\sigma$}} &\centering{\textbf{$time(s)$}} &\\
\hline
E6H0   & 6  & 7   &  & \centering{0.29446} &  & \textbf{0.2946670216} & \textbf{0.2946670216} & \textbf{0.2946670216} & 10/10 & 0.00 & 0.09   &  \\
E6H0   & 6  & 8   &  & \centering{0.28084} &  & \textbf{0.2810468468} & \textbf{0.2810468468} & \textbf{0.2810468468} & 10/10 & 0.00 & 0.09   &  \\
E6H0   & 6  & 9   &  & \centering{0.27273} &  & \textbf{0.2729182718} & \textbf{0.2729182718} & \textbf{0.2729182718} & 10/10 & 0.00 & 0.17   &  \\
E6H0   & 6  & 10  &  & \centering{0.26200} &  & \textbf{0.2621819240} & \textbf{0.2621819240} & \textbf{0.2621819240} & 10/10 & 0.00 & 0.17   &  \\
E6H0   & 6  & 11  &  & \centering{0.25017} &  & \textbf{0.2543330951} & \textbf{0.2543330951} & \textbf{0.2543330951} & 10/10 & 0.00 & 0.15   &  \\
E6H0   & 6  & 12  &  & \centering{0.24998} &  & \textbf{0.2500000000} & \textbf{0.2500000000} & \textbf{0.2500000000} & 10/10 & 0.00 & 0.14   &  \\
E6H0   & 6  & 13  &  & \centering{0.22667} &  & \textbf{0.2269506117} & \textbf{0.2269506117} & \textbf{0.2269506117} & 10/10 & 0.00 & 0.26   &  \\
E6H0   & 6  & 14  &  & \centering{0.22001} &  & \textbf{0.2201214487} & \textbf{0.2201214487} & \textbf{0.2201214487} & 10/10 & 0.00 & 1.25   &  \\
E6H0   & 6  & 15  &  & \centering{0.21153} &  & \textbf{0.2124800251} & \textbf{0.2124800251} & \textbf{0.2124800251} & 10/10 & 0.00 & 0.27   &  \\
E6H0   & 6  & 16  &  & \centering{0.20535} &  & \textbf{0.2075604739} & \textbf{0.2075604739} & \textbf{0.2075604739} & 10/10 & 0.00 & 0.36   &  \\
E11H0  & 11 & 70  &  & \centering{0.36610} &  & \textbf{0.3729322173} & \textbf{0.3729322173} & \textbf{0.3729322173} & 10/10 & 0.00 & 22.77  &  \\
E12H0b & 12 & 100 &  & \centering{0.32867} &  & \textbf{0.3361308135} & \textbf{0.3361308135} & \textbf{0.3361308135} & 10/10 & 0.00 & 351.41 &  \\
\hline
\#Best &    &     &  &   \centering{0} &  & \centering{12} & \centering{12} & \centering{12} &       &       &       &  \\
\hline
\end{tabular}
\end{tiny}
\end{table}
\renewcommand{\baselinestretch}{1.0}\large\normalsize

\begin{table*}\centering
\renewcommand{\arraystretch}{1.0}
\caption{Comparisons of the proposed TSGO algorithm with two reference algorithms (i.e., BH$^*$ and MBH$^*$) on the equal circle packing problem for problem instances with a region that consists of a small number of edges, where the best results among the compared algorithms are indicated in bold in terms of $R_{best}$ and $R_{avg}$.} \label{table_comparison1}
\begin{tiny}
\begin{tabular}{p{1.2cm}p{0.4cm}p{0.01cm}p{1.8cm} p{1.8cm}p{1.8cm} p{0.01cm}p{1.8cm} p{1.8cm} p{1.8cm}p{0.01cm} }
\hline
         &    & & \multicolumn{3}{c}{$R_{best}$} & &\multicolumn{3}{c}{$R_{avg}$} &\\
\cline{4-6}\cline{8-10}
Container    & $p$        &  &  BH$^*$  &  MBH$^*$ &  TSGO      &  &  BH$^*$ &  MBH$^*$ & TSGO & \\
\hline
E4H0a   & 50  &  & \textbf{0.0713771039} & \textbf{0.0713771039} & \textbf{0.0713771039} &  & \textbf{0.0713771039} & \textbf{0.0713771039} & \textbf{0.0713771039} &  \\
E4H0b   & 50  &  & \textbf{0.2721671717} & \textbf{0.2721671717} & \textbf{0.2721671717} &  & \textbf{0.2721671717} & \textbf{0.2721671717} & \textbf{0.2721671717} &  \\
E6H0    & 50  &  & \textbf{0.1226935182} & \textbf{0.1226935182} & \textbf{0.1226935182} &  & \textbf{0.1226935182} & \textbf{0.1226935182} & \textbf{0.1226935182} &  \\
E9H2    & 50  &  & 0.3686157489          & 0.3656662863          & \textbf{0.3689738179} &  & 0.3634125564          & 0.3616720252          & \textbf{0.3689738179} &  \\
E11H0   & 50  &  & \textbf{0.4368343737} & \textbf{0.4368343737} & \textbf{0.4368343737} &  & 0.4357490250          & 0.4365780965          & \textbf{0.4368343737} &  \\
E12H0a  & 50  &  & \textbf{0.6225798476} & \textbf{0.6225798476} & \textbf{0.6225798476} &  & 0.6225642581          & 0.6225763213          & \textbf{0.6225798476} &  \\
E12H0b  & 50  &  & \textbf{0.4620855791} & \textbf{0.4620855791} & \textbf{0.4620855791} &  & 0.4619205142          & 0.4618933105          & \textbf{0.4620855791} &  \\
E12H2   & 50  &  & \textbf{0.4448612715} & \textbf{0.4448612715} & \textbf{0.4448612715} &  & 0.4448604939          & 0.4448608827          & \textbf{0.4448612715} &  \\
E18H0   & 50  &  & \textbf{0.7877071509} & \textbf{0.7877071509} & \textbf{0.7877071509} &  & \textbf{0.7877071509} & \textbf{0.7877071509} & \textbf{0.7877071509} &  \\
E20H1   & 50  &  & \textbf{0.4137431060} & \textbf{0.4137431060} & \textbf{0.4137431060} &  & 0.4137379287          & 0.4137363551          & \textbf{0.4137431060} &  \\
E20H2   & 50  &  & \textbf{0.6402224988} & \textbf{0.6402224988} & \textbf{0.6402224988} &  & \textbf{0.6402224988} & 0.6401546638          & \textbf{0.6402224988} &  \\
E23H1   & 50  &  & \textbf{0.4005485124} & \textbf{0.4005485124} & \textbf{0.4005485124} &  & \textbf{0.4005485124} & \textbf{0.4005485124} & \textbf{0.4005485124} &  \\
E27H1   & 50  &  & \textbf{0.4844567507} & \textbf{0.4840050329} & \textbf{0.4844567507} &  & 0.4843452715          & 0.4840050329          & \textbf{0.4844567507} &  \\
E44H1   & 50  &  & \textbf{0.5345742337} & 0.5345665428          & \textbf{0.5345742337} &  & 0.5344076418          & 0.5342859578          & \textbf{0.5345742337} &  \\
E59H1   & 50  &  & 0.0199320353          & 0.0198678555          & \textbf{0.0199947494} &  & 0.0198443767          & 0.0198127987          & \textbf{0.0199947494} &  \\
E4H0a   & 100 &  & \textbf{0.0514010718} & \textbf{0.0514010718} & \textbf{0.0514010718} &  & 0.0514010668          & \textbf{0.0514010718} & \textbf{0.0514010718} &  \\
E4H0b   & 100 &  & \textbf{0.1978021937} & \textbf{0.1978021937} & \textbf{0.1978021937} &  & \textbf{0.1978021937} & \textbf{0.1978021937} & \textbf{0.1978021937} &  \\
E6H0    & 100 &  & \textbf{0.0882825264} & \textbf{0.0882825264} & \textbf{0.0882825264} &  & \textbf{0.0882825264} & \textbf{0.0882825264} & \textbf{0.0882825264} &  \\
E9H2    & 100 &  & \textbf{0.2801620820} & \textbf{0.2801620820} & \textbf{0.2801620820} &  & \textbf{0.2801620820} & 0.2798302518          & \textbf{0.2801620820} &  \\
E11H0   & 100 &  & \textbf{0.3181468932} & \textbf{0.3181468932} & \textbf{0.3181468932} &  & 0.3181468837          & 0.3181468853          & \textbf{0.3181468932} &  \\
E12H0a  & 100 &  & \textbf{0.4517241097} & \textbf{0.4517241097} & \textbf{0.4517241097} &  & \textbf{0.4517241097} & \textbf{0.4517241097} & \textbf{0.4517241097} &  \\
E12H0b  & 100 &  & \textbf{0.3361308135} & \textbf{0.3361308135} & \textbf{0.3361308135} &  & \textbf{0.3361308135} & \textbf{0.3361308135} & \textbf{0.3361308135} &  \\
E12H2   & 100 &  & \textbf{0.3172501921} & \textbf{0.3172501921} & \textbf{0.3172501921} &  & 0.3172493170          & \textbf{0.3172501921} & \textbf{0.3172501921} &  \\
E18H0   & 100 &  & \textbf{0.5761713874} & \textbf{0.5761713874} & \textbf{0.5761713874} &  & \textbf{0.5761713874} & \textbf{0.5761713874} & \textbf{0.5761713874} &  \\
E20H1   & 100 &  & \textbf{0.2989547054} & \textbf{0.2989547054} & \textbf{0.2989547054} &  & 0.2989526815          & \textbf{0.2989547054} & \textbf{0.2989547054} &  \\
E20H2   & 100 &  & \textbf{0.4733830176} & 0.4733830079          & \textbf{0.4733830176} &  & 0.4733823640          & 0.4733827279          & \textbf{0.4733830176} &  \\
E23H1   & 100 &  & \textbf{0.2896330498} & \textbf{0.2896330498} & \textbf{0.2896330498} &  & 0.2896299959          & \textbf{0.2896324594} & 0.2896316545          &  \\
E27H1   & 100 &  & \textbf{0.3541756870} & \textbf{0.3541756870} & \textbf{0.3541756870} &  & \textbf{0.3541756870} & \textbf{0.3541756870} & \textbf{0.3541756870} &  \\
E44H1   & 100 &  & \textbf{0.3907223295} & \textbf{0.3907223295} & \textbf{0.3907223295} &  & \textbf{0.3907223295} & \textbf{0.3907223295} & \textbf{0.3907223295} &  \\
E59H1   & 100 &  & \textbf{0.0150349633} & \textbf{0.0150349633} & \textbf{0.0150349633} &  & 0.0150334227          & 0.0150349542          & \textbf{0.0150349588} &  \\
E4H0a   & 150 &  & \textbf{0.0421454577} & \textbf{0.0421454577} & \textbf{0.0421454577} &  & 0.0421453491          & 0.0421453491          & \textbf{0.0421454454} &  \\
E4H0b   & 150 &  & \textbf{0.1636890708} & 0.1636890405          & \textbf{0.1636890708} &  & 0.1636890465          & 0.1636890405          & \textbf{0.1636890708} &  \\
E6H0    & 150 &  & \textbf{0.0725873025} & \textbf{0.0725873025} & \textbf{0.0725873025} &  & 0.0725872801          & \textbf{0.0725873025} & \textbf{0.0725873025} &  \\
E9H2    & 150 &  & \textbf{0.2333481422} & \textbf{0.2333481422} & \textbf{0.2333481422} &  & 0.2333137881          & 0.2332790355          & \textbf{0.2333481422} &  \\
E11H0   & 150 &  & \textbf{0.2628105627} & \textbf{0.2628105627} & \textbf{0.2628105627} &  & 0.2628082575          & 0.2627997871          & \textbf{0.2628104872} &  \\
E12H0a  & 150 &  & \textbf{0.3730487798} & \textbf{0.3730487798} & \textbf{0.3730487798} &  & 0.3730279634          & 0.3730429057          & \textbf{0.3730484422} &  \\
E12H0b  & 150 &  & \textbf{0.2783387066} & \textbf{0.2783387066} & \textbf{0.2783387066} &  & 0.2783385307          & 0.2783384399          & \textbf{0.2783387066} &  \\
E12H2   & 150 &  & \textbf{0.2636360626} & \textbf{0.2636360626} & \textbf{0.2636360626} &  & 0.2636357700          & 0.2636228471          & \textbf{0.2636360626} &  \\
E18H0   & 150 &  & 0.4800548333          & 0.4800569423          & \textbf{0.4800573753} &  & 0.4800409940          & \textbf{0.4800544890} & 0.4800534988          &  \\
E20H1   & 150 &  & 0.2508850879          & 0.2508904348          & \textbf{0.2508905275} &  & 0.2508718500          & \textbf{0.2508899440} & 0.2508892317          &  \\
E20H2   & 150 &  & \textbf{0.3929380796} & \textbf{0.3929380796} & \textbf{0.3929380796} &  & 0.3929376328          & 0.3929379824          & \textbf{0.3929380796} &  \\
E23H1   & 150 &  & 0.2425848102          & 0.2425853555          & \textbf{0.2425858370} &  & 0.2425745664          & 0.2425802538          & \textbf{0.2425856154} &  \\
E27H1   & 150 &  & 0.2930400468          & \textbf{0.2930778541} & 0.2930756336          &  & 0.2929020446          & 0.2929866807          & \textbf{0.2930447613} &  \\
E44H1   & 150 &  & 0.3253694621          & \textbf{0.3253697234} & \textbf{0.3253697234} &  & 0.3253445564          & 0.3253641926          & \textbf{0.3253657613} &  \\
E59H1   & 150 &  & 0.0126020291          & 0.0126004282          & \textbf{0.0126021266} &  & 0.0125991512          & 0.0125977296          & \textbf{0.0125992932} &  \\
E4H0a   & 200 &  & 0.0366127743          & 0.0366127127          & \textbf{0.0366127989} &  & 0.0366045964          & 0.0365977879          & \textbf{0.0366083653} &  \\
E4H0b   & 200 &  & \textbf{0.1429882126} & \textbf{0.1429882126} & \textbf{0.1429882126} &  & 0.1429867411          & 0.1429880516          & \textbf{0.1429882126} &  \\
E6H0    & 200 &  & \textbf{0.0636129856} & \textbf{0.0636129856} & \textbf{0.0636129856} &  & 0.0636129785          & 0.0636129789          & \textbf{0.0636129856} &  \\
E9H2    & 200 &  & 0.2051749126          & 0.2051734607          & \textbf{0.2051780567} &  & 0.2051581717          & 0.2050920008          & \textbf{0.2051759042} &  \\
E11H0   & 200 &  & \textbf{0.2300479293} & 0.2300479269          & \textbf{0.2300479293} &  & 0.2300475234          & 0.2300465204          & \textbf{0.2300479293} &  \\
E12H0a  & 200 &  & 0.3303497674          & 0.3303513835          & \textbf{0.3303513838} &  & 0.3303312676          & 0.3303323760          & \textbf{0.3303506670} &  \\
E12H0b  & 200 &  & \textbf{0.2397150105} & 0.2397149858          & 0.2397150046          &  & 0.2397144563          & 0.2397146114          & \textbf{0.2397147555} &  \\
E12H2   & 200 &  & \textbf{0.2334882838} & 0.2334869495          & 0.2334870331          &  & \textbf{0.2334860073} & 0.2334026236          & 0.2334801734          &  \\
E18H0   & 200 &  & 0.4190034532          & 0.4191417996          & \textbf{0.4191458163} &  & 0.4188462333          & 0.4190123862          & \textbf{0.4190839637} &  \\
E20H1   & 200 &  & 0.2172626068          & 0.2172697806          & \textbf{0.2172791061} &  & 0.2172200262          & 0.2172338832          & \textbf{0.2172676640} &  \\
E20H2   & 200 &  & 0.3446792546          & 0.3447052767          & \textbf{0.3447056183} &  & 0.3446675117          & 0.3446854835          & \textbf{0.3446855508} &  \\
E23H1   & 200 &  & 0.2101428378          & \textbf{0.2102077787} & 0.2102037845          &  & 0.2101305311          & 0.2101637466          & \textbf{0.2101661510} &  \\
E27H1   & 200 &  & \textbf{0.2595897715} & \textbf{0.2595897715} & \textbf{0.2595897715} &  & 0.2595276865          & 0.2595432008          & \textbf{0.2595486576} &  \\
E44H1   & 200 &  & 0.2813674161          & 0.2814511662          & \textbf{0.2816788437} &  & 0.2812891330          & 0.2814024292          & \textbf{0.2815233736} &  \\
E59H1   & 200 &  & 0.0111175198          & 0.0111161222          & \textbf{0.0111178733} &  & 0.0111155130          & 0.0111136582          & \textbf{0.0111156356} &  \\
\hline
\#Best  &     &  & 43                    & 40                    & \textbf{56}           &  & 15                    & 19                     & \textbf{56}      &  \\
\textit{p-value} &     &  & 4.63E-04              & 9.16E-04              &                       &  & 1.67E-08              & 1.25E-07              &                       &  \\
\hline
\end{tabular}
\end{tiny}
\end{table*}

\begin{table*}\centering
\renewcommand{\arraystretch}{1.0}
\caption{Comparisons of the proposed TSGO algorithm with two reference algorithms (i.e., BH$^*$ and MBH$^*$) on the equal circle packing problem for problem instances with a complicated region consisting of a large number of edges, where the best results among the compared algorithms are indicated in bold in terms of $R_{best}$ and $R_{avg}$.} \label{table_comparison1_complex}
\begin{tiny}
\begin{tabular}{p{1.2cm}p{0.4cm}p{0.01cm}p{1.8cm} p{1.8cm}p{1.8cm} p{0.01cm}p{1.8cm} p{1.8cm} p{1.8cm}p{0.01cm} }
\hline
         &    & & \multicolumn{3}{c}{$R_{best}$} & &\multicolumn{3}{c}{$R_{avg}$} &\\
\cline{4-6}\cline{8-10}
Container    & $p$        &  &  BH$^*$  &  MBH$^*$ &  TSGO      &  &  BH$^*$ &  MBH$^*$ & TSGO & \\
\hline
E101H2  & 100 &  & \textbf{0.5131508522} & \textbf{0.5131508522} & \textbf{0.5131508522} &  & 0.5131462111          & 0.5131375011          & \textbf{0.5131508522} &  \\
E101H3  & 100 &  & 0.7236466019          & 0.7238754116          & \textbf{0.7238846144} &  & 0.7233832936          & 0.7237614097          & \textbf{0.7238775278} &  \\
E106H3  & 100 &  & 0.6039146219          & \textbf{0.6039167191} & \textbf{0.6039167191} &  & 0.6037780899          & 0.6038959925          & \textbf{0.6039165094} &  \\
E106H5  & 100 &  & 0.5322579826          & 0.5328329107          & \textbf{0.5332531100} &  & 0.5308443285          & 0.5323526847          & \textbf{0.5331975141} &  \\
E107H3  & 100 &  & \textbf{1.1937852368} & \textbf{1.1937852368} & \textbf{1.1937852368} &  & \textbf{1.1937852368} & \textbf{1.1937852368} & \textbf{1.1937852368} &  \\
E120H3  & 100 &  & 0.8773504855          & 0.8773504855          & \textbf{0.8773789961} &  & 0.8772410287          & 0.8772947645          & \textbf{0.8773533279} &  \\
E172H4  & 100 &  & \textbf{0.4001953141} & 0.3970925656          & \textbf{0.4001953141} &  & 0.4001810712          & 0.3954351442          & \textbf{0.4001905629} &  \\
E193H1  & 100 &  & 0.8960463038          & \textbf{0.8970410179} & 0.8969499567          &  & 0.8953314693          & 0.8964794893          & \textbf{0.8966769264} &  \\
E196H5  & 100 &  & 0.4727866982          & 0.4652219431          & \textbf{0.4756223520} &  & 0.4709323222          & 0.4595169179          & \textbf{0.4756214961} &  \\
E203H3  & 100 &  & 0.2954793213          & 0.2892741421          & \textbf{0.2964789245} &  & 0.2932913997          & 0.2862528842          & \textbf{0.2957012517} &  \\
E60H1   & 100 &  & 0.1389364842          & \textbf{0.1389366502} & 0.1389364842          &  & 0.1389127305          & 0.1389155735          & \textbf{0.1389224366} &  \\
E81H3   & 100 &  & 0.4596181712          & \textbf{0.4596304590} & \textbf{0.4596304590} &  & 0.4595817333          & \textbf{0.4596304590} & 0.4596270833          &  \\
E82H3   & 100 &  & \textbf{0.4239063629} & \textbf{0.4239299929} & \textbf{0.4239299929} &  & 0.4238394299          & 0.4239298328          & \textbf{0.4239299929} &  \\
E84H3   & 100 &  & \textbf{0.4320481250} & \textbf{0.4320481250} & \textbf{0.4320481250} &  & \textbf{0.4320481250} & \textbf{0.4320481250} & \textbf{0.4320481250} &  \\
E93H4   & 100 &  & 0.5936319873          & 0.5936386605          & \textbf{0.5936444451} &  & 0.5934827880          & \textbf{0.5936305684} & 0.5936093706          &  \\
E101H2  & 150 &  & 0.4246039525          & 0.4246094078          & \textbf{0.4246301877} &  & 0.4245003932          & 0.4244980399          & \textbf{0.4245517063} &  \\
E101H3  & 150 &  & 0.6022903843          & 0.6022967719          & \textbf{0.6022981927} &  & 0.6022903843          & 0.6022689505          & \textbf{0.6022931886} &  \\
E106H3  & 150 &  & 0.4994823977          & 0.4996271982          & \textbf{0.4996758005} &  & 0.4994604604          & 0.4992973653          & \textbf{0.4995933394} &  \\
E106H5  & 150 &  & 0.4443063213          & \textbf{0.4444220001} & 0.4444091109          &  & 0.4441837084          & 0.4442980251          & \textbf{0.4443309623} &  \\
E107H3  & 150 &  & 0.9872584433          & 0.9872061372          & \textbf{0.9872881106} &  & 0.9871599470          & 0.9871103817          & \textbf{0.9872064963} &  \\
E120H3  & 150 &  & \textbf{0.7306445570} & 0.7306445568          & \textbf{0.7306445570} &  & 0.7306208580          & 0.7303939124          & \textbf{0.7306432422} &  \\
E172H4  & 150 &  & 0.3331491441          & 0.3331649258          & \textbf{0.3332646035} &  & 0.3330541489          & 0.3331410341          & \textbf{0.3331663193} &  \\
E193H1  & 150 &  & 0.7408977499          & \textbf{0.7409843995} & 0.7409137632          &  & \textbf{0.7406671323} & 0.7405549328          & 0.7404794296          &  \\
E196H5  & 150 &  & 0.3973150431          & 0.3971790864          & \textbf{0.3997206613} &  & 0.3955168408          & 0.3959672532          & \textbf{0.3995365485} &  \\
E203H3  & 150 &  & 0.2501434557          & 0.2483444155          & \textbf{0.2512178335} &  & 0.2495857093          & 0.2473432994          & \textbf{0.2500754648} &  \\
E60H1   & 150 &  & 0.1151603319          & \textbf{0.1151655673} & 0.1151505241          &  & 0.1151416702          & \textbf{0.1151512904} & 0.1151269139          &  \\
E81H3   & 150 &  & 0.3801972735          & \textbf{0.3802772848} & 0.3802766660          &  & 0.3801341086          & 0.3802365212          & \textbf{0.3802496299} &  \\
E82H3   & 150 &  & 0.3490790239          & 0.3490480371          & \textbf{0.3490791105} &  & 0.3490590187          & 0.3490358746          & \textbf{0.3490722172} &  \\
E84H3   & 150 &  & 0.3579755298          & 0.3579715432          & \textbf{0.3579871420} &  & 0.3578167382          & 0.3578566666          & \textbf{0.3579385074} &  \\
E93H4   & 150 &  & 0.4948046173          & \textbf{0.4949205191} & 0.4948698733          &  & 0.4947339911          & \textbf{0.4948543190} & 0.4948247350          &  \\
E101H2  & 200 &  & 0.3707706478          & 0.3708004504          & \textbf{0.3708266880} &  & 0.3706979030          & 0.3706963168          & \textbf{0.3708085633} &  \\
E101H3  & 200 &  & 0.5271335274          & 0.5268009361          & \textbf{0.5271473970} &  & 0.5268565459          & 0.5265542452          & \textbf{0.5270959759} &  \\
E106H3  & 200 &  & 0.4384135384          & 0.4384076393          & \textbf{0.4384167531} &  & 0.4383338932          & 0.4382993999          & \textbf{0.4383871870} &  \\
E106H5  & 200 &  & 0.3909836192          & 0.3909721826          & \textbf{0.3910559351} &  & 0.3909032381          & 0.3908447285          & \textbf{0.3910204812} &  \\
E107H3  & 200 &  & 0.8600500229          & \textbf{0.8601675879} & 0.8601603145          &  & 0.8598515845          & \textbf{0.8600649578} & 0.8600083142          &  \\
E120H3  & 200 &  & 0.6352015578          & 0.6352180420          & \textbf{0.6352456666} &  & 0.6349999088          & 0.6351366779          & \textbf{0.6351595315} &  \\
E172H4  & 200 &  & 0.2937985244          & 0.2938115968          & \textbf{0.2938264340} &  & \textbf{0.2936556747} & 0.2935779286          & 0.2935365465          &  \\
E193H1  & 200 &  & 0.6495304786          & 0.6496346640          & \textbf{0.6497271655} &  & 0.6493698892          & 0.6494441862          & \textbf{0.6495499522} &  \\
E196H5  & 200 &  & 0.3528640237          & 0.3532782171          & \textbf{0.3532967101} &  & 0.3526929504          & 0.3529051796          & \textbf{0.3532569447} &  \\
E203H3  & 200 &  & 0.2218022505          & 0.2230280767          & \textbf{0.2233159947} &  & 0.2214404497          & 0.2218444831          & \textbf{0.2232194868} &  \\
E60H1   & 200 &  & 0.1002412092          & 0.1002605523          & \textbf{0.1002674230} &  & 0.1002173527          & 0.1002302451          & \textbf{0.1002478734} &  \\
E81H3   & 200 &  & 0.3321541887          & 0.3322738793          & \textbf{0.3322827528} &  & 0.3320943976          & 0.3321386094          & \textbf{0.3321969424} &  \\
E82H3   & 200 &  & \textbf{0.3045767624} & 0.3045651214          & 0.3045495631          &  & \textbf{0.3044992578} & 0.3044837034          & 0.3044844759          &  \\
E84H3   & 200 &  & 0.3144020791          & 0.3143943149          & \textbf{0.3144135487} &  & 0.3143300047          & 0.3142707651          & \textbf{0.3143854408} &  \\
E93H4   & 200 &  & 0.4321912739          & 0.4321330929          & \textbf{0.4323742595} &  & 0.4320001139          & 0.4318109785          & \textbf{0.4322075414} &  \\
\hline
\#Best  &     &  & 7                     & 14                    & \textbf{36}           &  & 5                     & 7                     & \textbf{37}           &  \\
\textit{p-value} &     &  & 2.42E-07     & 2.42E-04              &                       &  & 1.01E-06              & 5.30E-06              &                       &  \\
\hline
\end{tabular}
\end{tiny}
\end{table*}

To further evaluate TSGO's performance on more complicated instances with holes, we carried out another experiment based on two sets of additional benchmark instances. The first set consists of 60 instances where the regions to be packed are relatively simple with a small number of edges and holes containing up to $p=200$ dispersion points. In particular, the number $|H|$ of holes in the region is at most two for each instance. The second set consists of 45 instances containing up to 200 dispersion points, where the regions to be packed are much more complex with a large number of edges and holes. In this experiment, we created two TSGO variants BH$^*$ and MBH$^*$ by respectively replacing the tabu search method of TSGO with the popular basin-hoping (BH) algorithm \citep{wales1997global} and the monotonic basin-hopping (MBH) algorithm \citep{leary2000global}, while keeping other TSGO components unchanged. Thus, BH$^*$ and MBH$^*$ employ the same optimization model as the TSGO algorithm, allowing us to make a fair comparison. For the subroutine BH of BH$^*$, the maximum number of iterations was set to 1000 for a single run and the temperature parameter $t$ was set to 0.1, and for the subroutine MBH of MBH$^*$ the search depth $\theta_{max}$ was set to 100 to ensure an intensified search.

The computational results of the BH$^*$, MBH$^*$, and TSGO algorithms are summarized in Tables \ref{table_comparison1} and \ref{table_comparison1_complex} respectively for the first and second set of instances, including the best objective value ($R_{best}$) over 10 independent runs and the average objective value ($R_{avg}$). The row ``\#Best"of each table indicates the numbers of instances for which the corresponding algorithm yields the best result among the compared algorithms in terms of $R_{best}$ and $R_{avg}$. The \textit{p-values} from the Wilcoxon signed-rank test are provided in the last row to show the statistical difference between the TSGO algorithm and the BH$^*$ and MBH$^*$ algorithms in terms of $R_{best}$ and $R_{avg}$, where a \textit{p-value} less than 0.05 means that there exists a significant difference between two groups of compared results. 

Table \ref{table_comparison1} shows that the TSGO algorithm significantly outperforms BH$^*$ and MBH$^*$  on both performance indicators. In terms of $R_{best}$, BH$^*$, MBH$^*$ and TSGO respectively obtained the best result for 43, 40 and 56 out of the 60 instances. In terms of $R_{avg}$, the superiority of the TSGO algorithm over the BH$^*$ and MBH$^*$ algorithms is even more pronounced, where the number of instances for which the TSGO algorithm obtained the best result (56) is much larger than the numbers of instances (15 and 19) for which the BH$^*$ and MBH$^*$ algorithms obtained the best result. The small \textit{p-values} ($\le 0.05$) additionally confirm that the differences between the TSGO algorithm and the BH$^*$ and MBH$^*$ algorithms are statistically significant for both $R_{best}$ and $R_{avg}$.

Table \ref{table_comparison1_complex} further shows that for the more complicated instances the TSGO algorithm performs much better than BH$^*$ and MBH$^*$. TSGO's superiority for these instances is more conspicuous compared to the previous instances in terms of both $R_{best}$ and $R_{avg}$. In terms of $R_{best}$, BH$^*$, MBH$^*$ and TSGO respectively obtained the best result for 7, 14 and 36 out of 45 instances. In terms of $R_{avg}$, the number of instances for which the TSGO algorithm obtained the best result is 37 which is much larger than that (5 and 7) of the BH$^*$ and MBH$^*$ algorithms. This experiment indicates that the proposed TSGO algorithm is particularly efficient for the equal circle packing problem in a complicated region composed of an irregular container and a number of irregular holes.

\begin{figure*}[h]
\centering
\subfigure[E60H1 for $p=200$]{\includegraphics[width=1.6in]{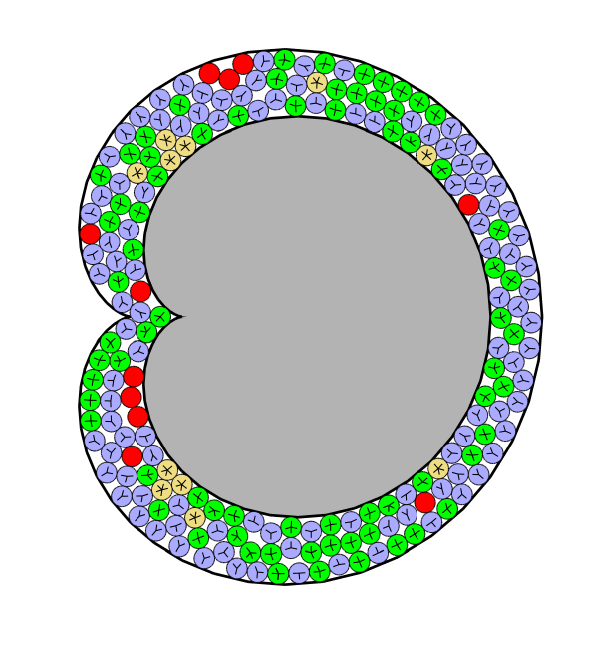}}
\subfigure[E101H2 for $p=200$]{\includegraphics[width=1.6in]{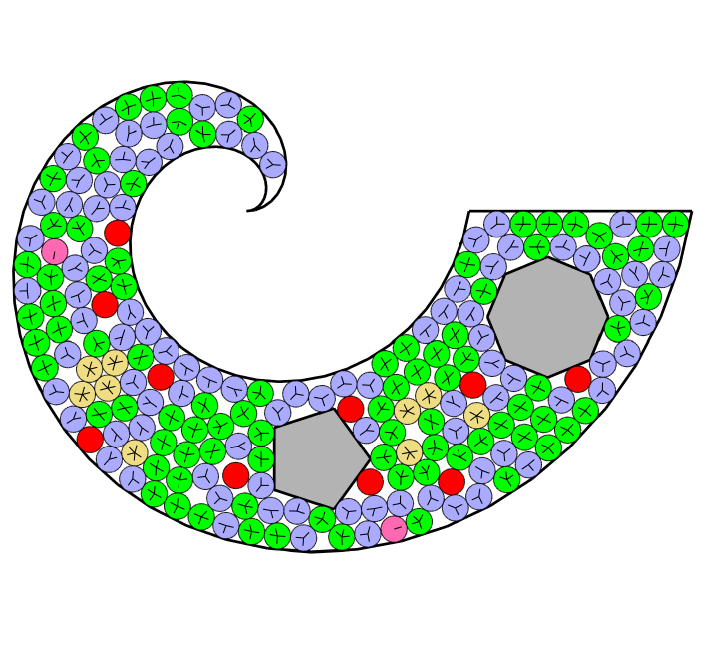}}
\subfigure[E193H1 for $p=200$]{\includegraphics[width=1.6in]{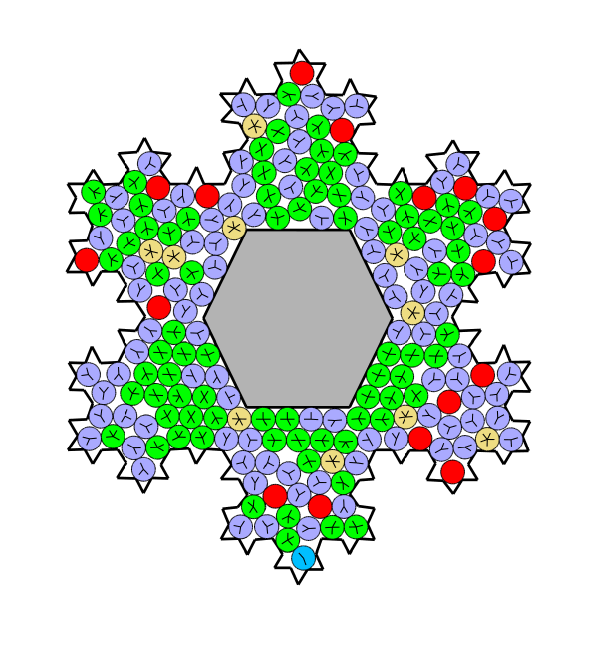}}
\subfigure[E106H5 for $p=200$]{\includegraphics[width=1.6in]{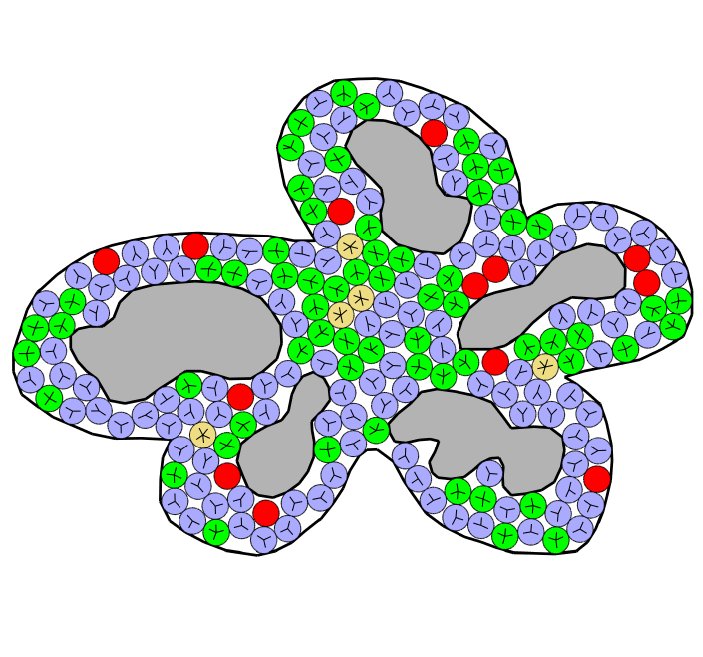}}
\subfigure[E120H3 for $p=200$]{\includegraphics[width=1.6in]{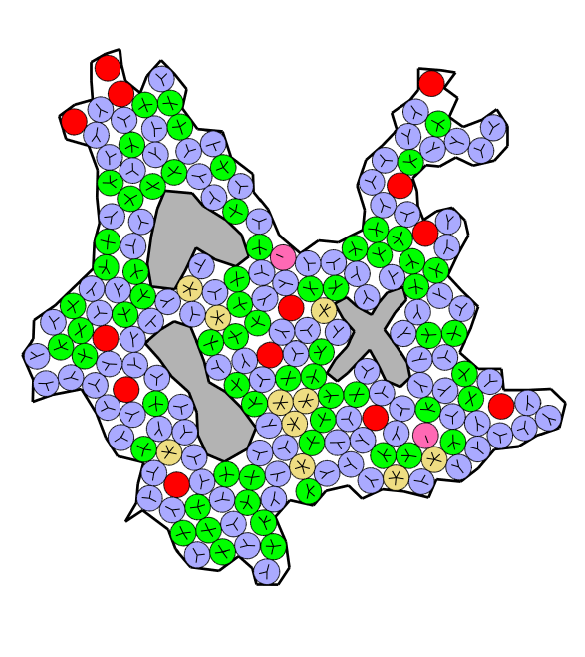}}
\subfigure[E9H2  for $p=200$]{\includegraphics[width=1.6in]{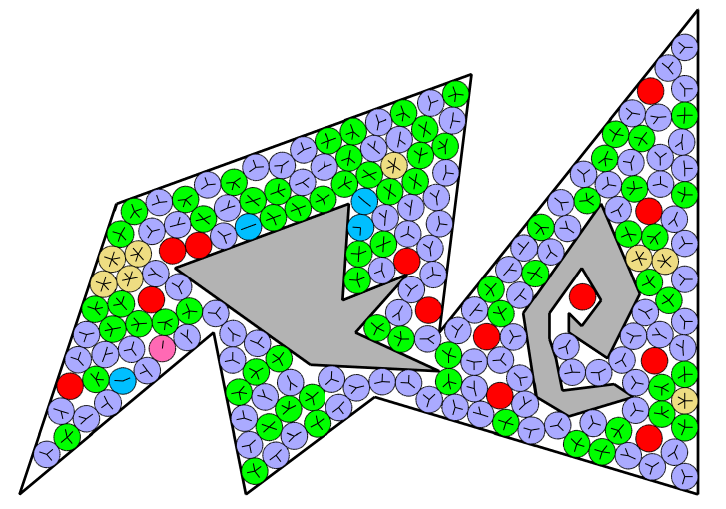}}
\subfigure[E84H3 for $p=200$]{\includegraphics[width=1.61in]{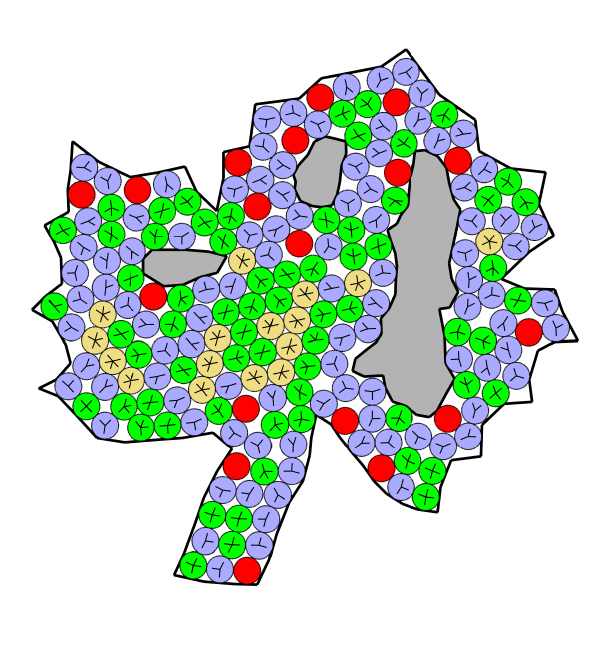}}
\subfigure[E59H1 for $p=150$]{\includegraphics[width=1.6in]{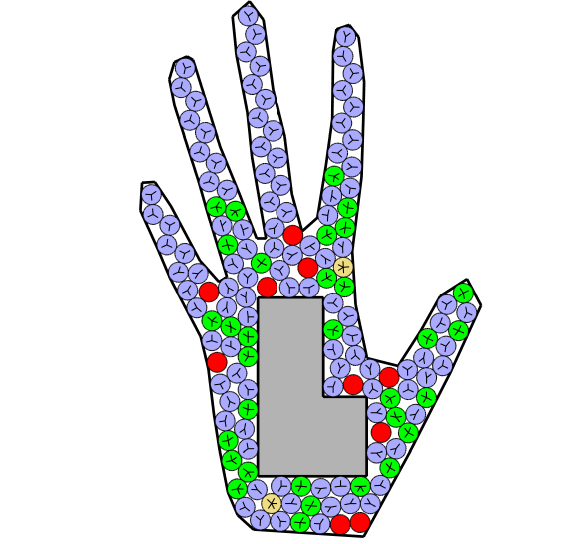}}
\subfigure[E203H3 for $p=200$]{\includegraphics[width=1.6in]{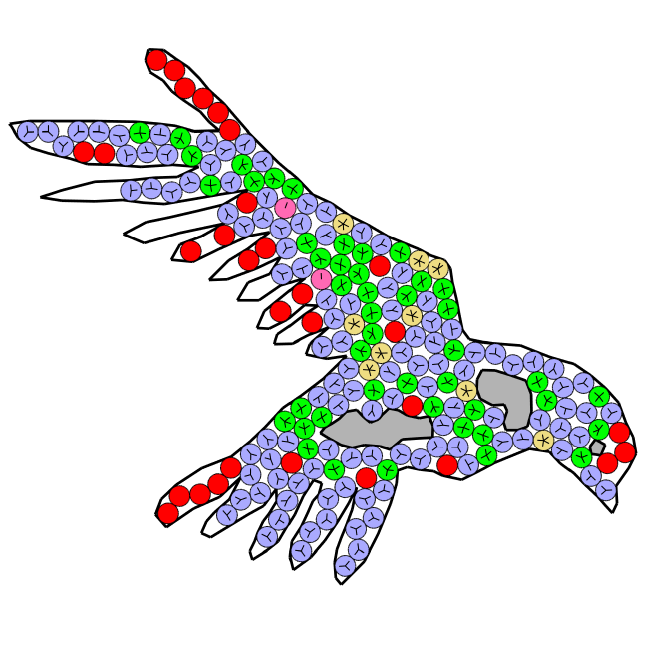}}
\caption{Best configurations found in this work for 9 representative instances, where the circles are colored according to their number of neighbors.}
\label{packing_improved_large}
\end{figure*}

Fig. \ref{packing_improved_large} gives graphical representations to provide an intuitive impression of the best configurations found for some representative instances. We see that the proposed model and algorithm are capable of effectively tackling highly complicated instances, including those containing a number of irregular holes and those composed of a large number of edges. The figures also show that the neighborhood operation of the tabu search, which moves high-energy dispersion points from their current positions to low-energy vacant sites, is very appropriate compared to the popular random displacement operation. Taking E9H2 with $p=200$ as an example, it is very difficult for the sightless random displacement operation to reach the current best configuration (i.e., the subfigure (f)) when the dent of the right-hand polygonal hole does not surround a dispersion point in the initial solution.

\subsection{Computational results and comparison on the point arrangement problem}
\label{results_withoutBC}

\begin{table*}\centering
\renewcommand{\arraystretch}{1.0}
\caption{Comparisons of the proposed TSGO algorithm with two reference algorithms (i.e., BH$^*$ and MBH$^*$) on the point arrangement problem for problem instances with a region that consists of a small number of edges, where the best results among three compared algorithms are indicated in bold both in terms of $D_{best}$ and $D_{avg}$.} \label{table_comparison2}
\begin{tiny}
\begin{tabular}{p{1.2cm}p{0.4cm}p{0.01cm}p{1.8cm} p{1.8cm}p{1.8cm} p{0.01cm}p{1.8cm} p{1.8cm} p{1.8cm}p{0.01cm} }
\hline
         &    & & \multicolumn{3}{c}{$D_{best}$} & &\multicolumn{3}{c}{$D_{avg}$} &\\
\cline{4-6}\cline{8-10}
Container    & $p$        &  &  BH$^*$  &  MBH$^*$ &  TSGO      &  &  BH$^*$ &  MBH$^*$ & TSGO & \\
\hline
E4H0a   & 50  &  & \textbf{0.1664988509} & \textbf{0.1664988509} & \textbf{0.1664988509} &  & \textbf{0.1664988509} & \textbf{0.1664988509} & \textbf{0.1664988509} &  \\
E4H0b   & 50  &  & \textbf{0.7232436230} & \textbf{0.7232436230} & \textbf{0.7232436230} &  & \textbf{0.7232436230} & \textbf{0.7232436230} & \textbf{0.7232436230} &  \\
E6H0    & 50  &  & 0.2933779419          & 0.2933779419          & \textbf{0.2933786207} &  & 0.2933779419          & 0.2933779419          & \textbf{0.2933780098} &  \\
E9H2    & 50  &  & \textbf{1.1252199453} & \textbf{1.1252199453} & \textbf{1.1252199453} &  & 1.1250991089          & \textbf{1.1252199453} & \textbf{1.1252199453} &  \\
E11H0   & 50  &  & \textbf{1.1330657353} & \textbf{1.1330657353} & \textbf{1.1330657353} &  & 1.1330162295          & \textbf{1.1330657353} & \textbf{1.1330657353} &  \\
E12H0a  & 50  &  & \textbf{1.6364303754} & \textbf{1.6364303754} & \textbf{1.6364303754} &  & \textbf{1.6364303754} & \textbf{1.6364303754} & \textbf{1.6364303754} &  \\
E12H0b  & 50  &  & \textbf{1.1512460802} & \textbf{1.1512460802} & \textbf{1.1512460802} &  & \textbf{1.1512460802} & \textbf{1.1512460802} & \textbf{1.1512460802} &  \\
E12H2   & 50  &  & \textbf{1.2090410970} & \textbf{1.2090410970} & \textbf{1.2090410970} &  & 1.2089809516          & \textbf{1.2090410970} & \textbf{1.2090410970} &  \\
E18H0   & 50  &  & \textbf{2.2664436736} & \textbf{2.2664436736} & \textbf{2.2664436736} &  & 2.2664388856          & \textbf{2.2664436736} & \textbf{2.2664436736} &  \\
E20H1   & 50  &  & \textbf{1.1494074146} & \textbf{1.1494074146} & \textbf{1.1494074146} &  & 1.1487617666          & 1.1483431535          & \textbf{1.1492486621} &  \\
E20H2   & 50  &  & 1.7642062858          & \textbf{1.7642062920} & 1.7642062877          &  & 1.7640056517          & \textbf{1.7642062866} & \textbf{1.7642062866} &  \\
E23H1   & 50  &  & \textbf{1.0945028103} & \textbf{1.0945028103} & \textbf{1.0945028103} &  & 1.0944768737          & 1.0943600422          & \textbf{1.0944875530} &  \\
E27H1   & 50  &  & \textbf{1.3363883975} & \textbf{1.3363883975} & \textbf{1.3363883975} &  & 1.3353178683          & \textbf{1.3363490513} & \textbf{1.3363490513} &  \\
E44H1   & 50  &  & \textbf{1.4231184848} & \textbf{1.4231184848} & \textbf{1.4231184848} &  & 1.4229226152          & 1.4230690380          & \textbf{1.4231165211} &  \\
E59H1   & 50  &  & 0.0625113933          & \textbf{0.0625365968} & \textbf{0.0625365968} &  & 0.0624453409          & 0.0624835536          & \textbf{0.0625045164} &  \\
E4H0a   & 100 &  & 0.1145681001          & \textbf{0.1145682248} & \textbf{0.1145682248} &  & 0.1145662184          & \textbf{0.1145671956} & 0.1145670435          &  \\
E4H0b   & 100 &  & \textbf{0.4826373801} & \textbf{0.4826373801} & \textbf{0.4826373801} &  & 0.4826356082          & 0.4826179302          & \textbf{0.4826359562} &  \\
E6H0    & 100 &  & \textbf{0.1998574337} & 0.1998574216          & 0.1998574281          &  & 0.1998560311          & 0.1997843134          & \textbf{0.1998568131} &  \\
E9H2    & 100 &  & 0.7494661035          & 0.7511238475          & \textbf{0.7511278062} &  & 0.7478403599          & 0.7502858249          & \textbf{0.7505626930} &  \\
E11H0   & 100 &  & 0.7655396868          & 0.7656897992          & \textbf{0.7656953424} &  & 0.7652491178          & 0.7655114064          & \textbf{0.7656381539} &  \\
E12H0a  & 100 &  & \textbf{1.1078034265} & \textbf{1.1078034265} & \textbf{1.1078034265} &  & \textbf{1.1077982675} & 1.1077723402          & 1.1077925638          &  \\
E12H0b  & 100 &  & 0.7755043144          & 0.7755043105          & \textbf{0.7755043306} &  & 0.7754464078          & \textbf{0.7754653422} & 0.7754471255          &  \\
E12H2   & 100 &  & \textbf{0.8146928310} & 0.8146928235          & \textbf{0.8146928310} &  & 0.8143454056          & \textbf{0.8146115427} & 0.8145986367          &  \\
E18H0   & 100 &  & 1.4913313807          & \textbf{1.4935027529} & 1.4928956213          &  & 1.4902461911          & \textbf{1.4928152147} & 1.4920205838          &  \\
E20H1   & 100 &  & 0.7639769482          & 0.7638994850          & \textbf{0.7640895485} &  & 0.7636734678          & 0.7636264998          & \textbf{0.7638340934} &  \\
E20H2   & 100 &  & 1.1842551290          & 1.1842672680          & \textbf{1.1842695666} &  & 1.1829491272          & 1.1842389685          & \textbf{1.1842567749} &  \\
E23H1   & 100 &  & 0.7209507423          & 0.7211250210          & \textbf{0.7211499167} &  & 0.7196800678          & 0.7205693410          & \textbf{0.7207608891} &  \\
E27H1   & 100 &  & 0.9087618011          & \textbf{0.9088386862} & 0.9088386700          &  & 0.9063990417          & 0.9077689522          & \textbf{0.9087282509} &  \\
E44H1   & 100 &  & 0.9655365849          & 0.9655319844          & \textbf{0.9656074713} &  & 0.9633779411          & 0.9643480822          & \textbf{0.9649732229} &  \\
E59H1   & 100 &  & 0.0424090996          & \textbf{0.0424419006} & 0.0424207994          &  & 0.0423971176          & 0.0423749032          & \textbf{0.0423973095} &  \\
\hline
\#Best  &     &  & 16                    & 20                    & \textbf{25}           &  & 5                     & 14                    & \textbf{25}           &  \\
\textit{p-value} &     &  & 9.87E-04     & 1.40E-01              &                       &  & 2.10E-05              & 4.55E-03              &                       &  \\
\hline
\end{tabular}
\end{tiny}
\end{table*}

\begin{table*}\centering
\renewcommand{\arraystretch}{1.0}
\caption{Comparisons of the proposed TSGO algorithm with two reference algorithms (i.e., BH$^*$ and MBH$^*$) on the point arrangement problem for problem instances with a complicated region that consists of a large number of edges, where the best results among three compared algorithms are indicated in bold both in terms of $D_{best}$ and $D_{avg}$.} \label{table_comparison2_complex}
\begin{tiny}
\begin{tabular}{p{1.2cm}p{0.4cm}p{0.01cm}p{1.8cm} p{1.8cm}p{1.8cm} p{0.01cm}p{1.8cm} p{1.8cm} p{1.8cm}p{0.01cm} }
\hline
         &    & & \multicolumn{3}{c}{$D_{best}$} & &\multicolumn{3}{c}{$D_{avg}$} &\\
\cline{4-6}\cline{8-10}
Container    & $p$        &  &  BH$^*$  &  MBH$^*$ &  TSGO      &  &  BH$^*$ &  MBH$^*$ & TSGO & \\
\hline
E101H2  & 50  &  & 1.9294792183          & \textbf{1.9296611433} & \textbf{1.9296611433} &  & 1.7373305150 & \textbf{1.9293080575} & 1.9290832147          &  \\
E101H3  & 50  &  & \textbf{2.7597012147} & \textbf{2.7597012147} & \textbf{2.7597012147} &  & 2.7579227475 & 2.7528616370          & \textbf{2.7583284264} &  \\
E106H3  & 50  &  & 2.2760303651          & \textbf{2.2785596163} & \textbf{2.2785596163} &  & 2.2721291472 & 2.2753498323          & \textbf{2.2757095844} &  \\
E106H5  & 50  &  & 2.1492581495          & \textbf{2.1503429086} & \textbf{2.1503429086} &  & 2.1474806367 & 2.1495254107          & \textbf{2.1495595896} &  \\
E107H3  & 50  &  & \textbf{4.3776888646} & \textbf{4.3776888646} & \textbf{4.3776888646} &  & 4.3744222239 & 4.3763245600          & \textbf{4.3764749744} &  \\
E120H3  & 50  &  & \textbf{3.3713160251} & \textbf{3.3713160251} & 3.3711157887          &  & 3.3694128674 & \textbf{3.3707797719} & 3.3707419523          &  \\
E172H4  & 50  &  & \textbf{1.6086033153} & 1.6083403888          & 1.6083403888          &  & 1.6073993775 & \textbf{1.6080920152} & 1.6079020027          &  \\
E193H1  & 50  &  & 3.4885688530          & 3.4930481925          & \textbf{3.4931672197} &  & 3.4809669506 & 3.4828879612          & \textbf{3.4849529684} &  \\
E196H5  & 50  &  & 2.1462437211          & 2.1458081483          & \textbf{2.1525660024} &  & 2.1393500463 & 2.1359675846          & \textbf{2.1420151416} &  \\
E203H3  & 50  &  & 1.2947250241          & \textbf{1.2963445060} & \textbf{1.2963445060} &  & 1.2907711108 & 1.2944729939          & \textbf{1.2950615485} &  \\
E60H1   & 50  &  & \textbf{0.5509805177} & \textbf{0.5509805177} & \textbf{0.5509805177} &  & 0.5509455932 & 0.5508626165          & \textbf{0.5509709816} &  \\
E81H3   & 50  &  & \textbf{1.6902021702} & \textbf{1.6902021702} & \textbf{1.6902021702} &  & 1.6895559247 & 1.6900255476          & \textbf{1.6900946298} &  \\
E82H3   & 50  &  & 1.5476779639          & 1.5477194735          & \textbf{1.5477286761} &  & 1.5460645514 & 1.5455161767          & \textbf{1.5464949834} &  \\
E84H3   & 50  &  & \textbf{1.6205855793} & \textbf{1.6205855793} & \textbf{1.6205855793} &  & 1.6202744664 & 1.6202793679          & \textbf{1.6205050592} &  \\
E93H4   & 50  &  & 2.4057626982          & 2.4087452204          & \textbf{2.4096068833} &  & 2.4038957388 & 2.4047430949          & \textbf{2.4057504187} &  \\
E101H2  & 100 &  & \textbf{1.2840210583} & 1.2839958467          & 1.2833642210          &  & 1.2816485778 & \textbf{1.2822682646} & 1.2822230393          &  \\
E101H3  & 100 &  & 1.8474917504          & 1.8517629108          & \textbf{1.8527605603} &  & 1.8454387326 & 1.8493839200          & \textbf{1.8496615212} &  \\
E106H3  & 100 &  & 1.5368444673          & 1.5389609812          & \textbf{1.5397160821} &  & 1.5336969076 & 1.5373017679          & \textbf{1.5382859065} &  \\
E106H5  & 100 &  & 1.4471249445          & \textbf{1.4471430022} & \textbf{1.4471430022} &  & 1.4448386089 & 1.4449916914          & \textbf{1.4460133622} &  \\
E107H3  & 100 &  & 2.9391747678          & 2.9422207363          & \textbf{2.9423829506} &  & 2.9367746053 & 2.9394661810          & \textbf{2.9397676542} &  \\
E120H3  & 100 &  & 2.2566551682          & \textbf{2.2585254097} & 2.2581001199          &  & 2.2527411251 & \textbf{2.2548931409} & 2.2535581665          &  \\
E172H4  & 100 &  & 1.0758381251          & \textbf{1.0767384436} & 1.0759548576          &  & 1.0727718180 & \textbf{1.0749063077} & 1.0736695957          &  \\
E193H1  & 100 &  & 2.3254402936          & \textbf{2.3343111553} & \textbf{2.3343111553} &  & 2.3205631956 & \textbf{2.3295235036} & 2.3287152398          &  \\
E196H5  & 100 &  & 1.4098783673          & 1.4138861609          & \textbf{1.4209052194} &  & 1.4065769102 & 1.4103438073          & \textbf{1.4128738037} &  \\
E203H3  & 100 &  & 0.8522060778          & \textbf{0.8564731696} & 0.8557972355          &  & 0.8507183234 & 0.8539301996          & \textbf{0.8544934583} &  \\
E60H1   & 100 &  & 0.3582299427          & 0.3582466332          & \textbf{0.3582993202} &  & 0.3577676480 & 0.3578719449          & \textbf{0.3580180655} &  \\
E81H3   & 100 &  & 1.1482805892          & 1.1485277749          & \textbf{1.1485321215} &  & 1.1472277187 & 1.1471834164          & \textbf{1.1479214956} &  \\
E82H3   & 100 &  & 1.0473915944          & \textbf{1.0474433895} & 1.0474228732          &  & 1.0469480201 & 1.0470108082          & \textbf{1.0471543246} &  \\
E84H3   & 100 &  & 1.0869889031          & \textbf{1.0872842260} & 1.0871845694          &  & 1.0853519460 & 1.0860224488          & \textbf{1.0863459158} &  \\
E93H4   & 100 &  & 1.5934017314          & \textbf{1.5939701632} & 1.5937588256          &  & 1.5916298296 & 1.5925545584          & \textbf{1.5930971364} &  \\
E101H2  & 150 &  & 1.0166314840          & 1.0155772115          & \textbf{1.0175701224} &  & 1.0127636545 & 1.0143834682          & \textbf{1.0163457379} &  \\
E101H3  & 150 &  & 1.4749399669          & \textbf{1.4759039711} & 1.4758619447          &  & 1.4733512716 & 1.4737548666          & \textbf{1.4748651219} &  \\
E106H3  & 150 &  & 1.2301928882          & 1.2307360530          & \textbf{1.2316898249} &  & 1.2272249098 & 1.2276956567          & \textbf{1.2296099299} &  \\
E106H5  & 150 &  & 1.1379054440          & 1.1372298361          & \textbf{1.1389649453} &  & 1.1370513922 & 1.1364212055          & \textbf{1.1373346074} &  \\
E107H3  & 150 &  & 2.3545859328          & \textbf{2.3577371747} & 2.3574562564          &  & 2.3507354789 & 2.3531119086          & \textbf{2.3532942920} &  \\
E120H3  & 150 &  & 1.7957563567          & 1.7975153377          & \textbf{1.8007023797} &  & 1.7927535594 & 1.7955887570          & \textbf{1.7959811597} &  \\
E172H4  & 150 &  & 0.8524978245          & \textbf{0.8525962904} & 0.8516701412          &  & 0.8508764111 & 0.8509029887          & \textbf{0.8510752628} &  \\
E193H1  & 150 &  & 1.8619806126          & 1.8620687888          & \textbf{1.8630694963} &  & 1.8540394316 & 1.8561542253          & \textbf{1.8562745973} &  \\
E196H5  & 150 &  & 1.1085450824          & 1.1096023844          & \textbf{1.1096389040} &  & 1.1015314163 & 1.1007478968          & \textbf{1.1045592005} &  \\
E203H3  & 150 &  & 0.6798294129          & 0.6801848429          & \textbf{0.6828531990} &  & 0.6776296256 & 0.6786975696          & \textbf{0.6802848226} &  \\
E60H1   & 150 &  & 0.2827053641          & 0.2828959343          & \textbf{0.2829748122} &  & 0.2824164724 & 0.2826403084          & \textbf{0.2827246206} &  \\
E81H3   & 150 &  & 0.9151500861          & 0.9155216189          & \textbf{0.9163973051} &  & 0.9131879749 & 0.9141138269          & \textbf{0.9149487359} &  \\
E82H3   & 150 &  & 0.8314206381          & \textbf{0.8324378837} & 0.8321274553          &  & 0.8310173646 & 0.8314429821          & \textbf{0.8314888335} &  \\
E84H3   & 150 &  & \textbf{0.8630536076} & 0.8629312866          & 0.8629057017          &  & 0.8620981076 & 0.8618367930          & \textbf{0.8620988622} &  \\
E93H4   & 150 &  & \textbf{1.2604775170} & 1.2601949020          & 1.2604568525          &  & 1.2587040149 & 1.2588156611          & \textbf{1.2595900730} &  \\
\hline
\#Best  &     &  & 10                    & 22                    & 30                    &  & 0            & 7                     & 38                    &  \\
\textit{p-value} &     &  & 3.10E-06     & 6.44E-02              &                       &  & 5.18E-09     & 1.58E-05              &                       &  \\
\hline
\end{tabular}
\end{tiny}
\end{table*}

\begin{figure*}[h]
\centering
\centering
\subfigure[E60H1 for $p=100$]{\includegraphics[width=1.6in]{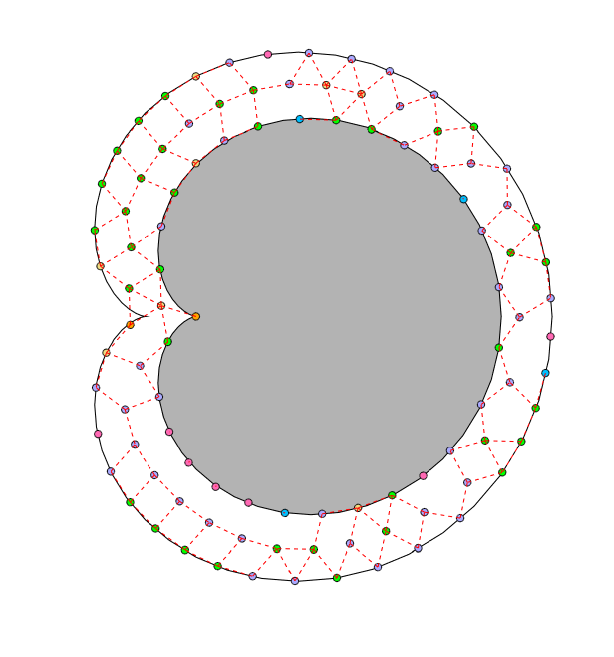}}
\subfigure[E101H2 for $p=100$]{\includegraphics[width=1.6in]{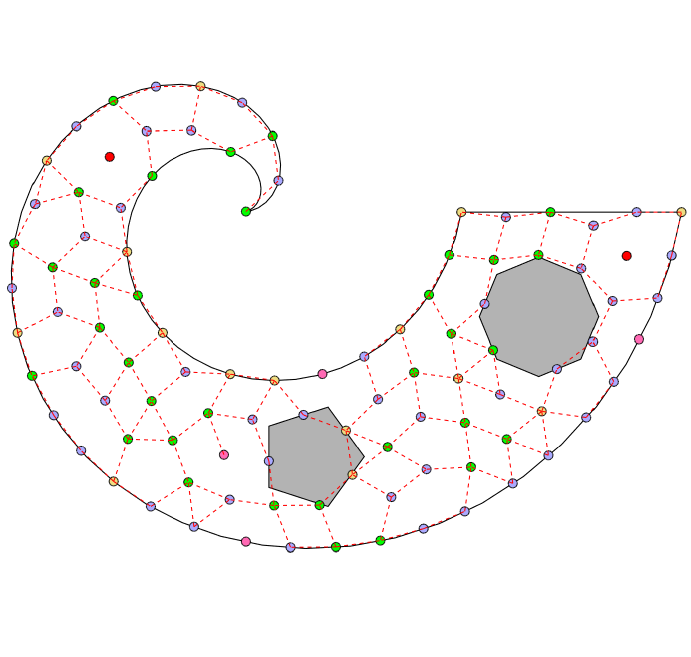}}
\subfigure[E44H1 for $p=100$]{\includegraphics[width=1.6in]{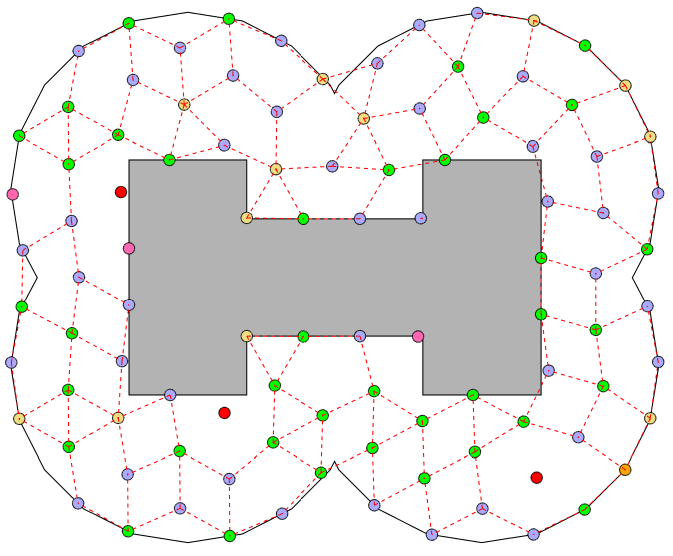}}
\subfigure[E84H3 for $p=100$]{\includegraphics[width=1.6in]{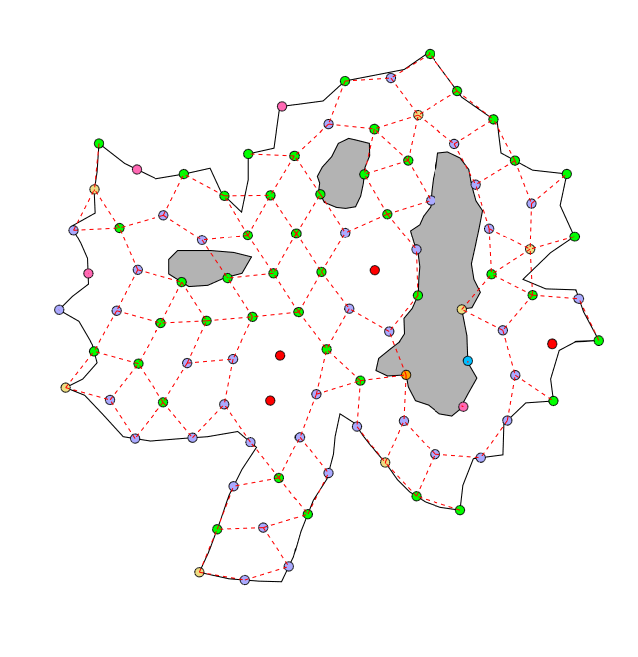}}
\subfigure[E93H4 for $p=100$]{\includegraphics[width=1.6in]{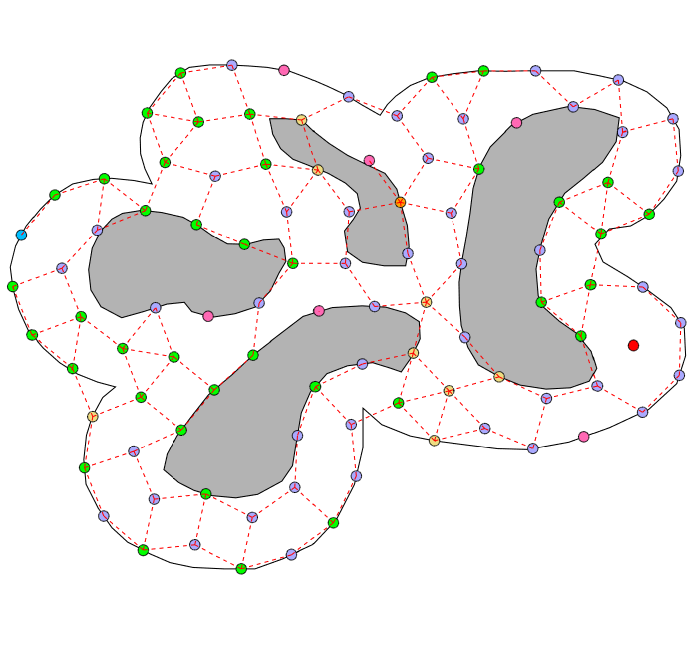}}
\subfigure[E107H3 for $p=100$]{\includegraphics[width=1.6in]{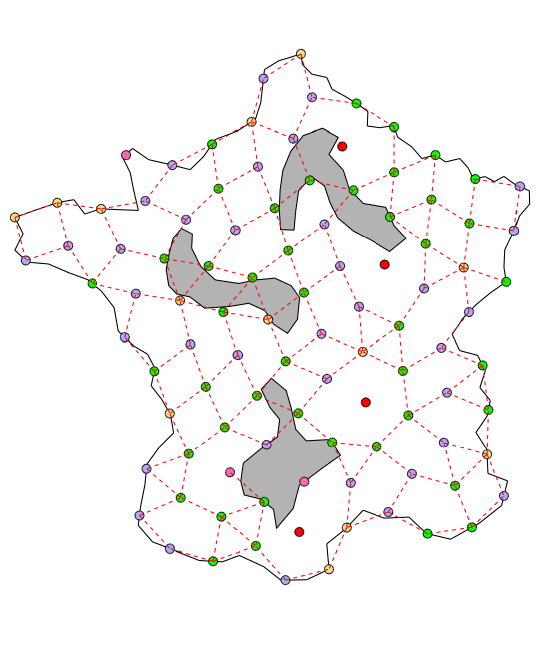}}
\subfigure[E120H3 for $p=150$]{\includegraphics[width=1.6in]{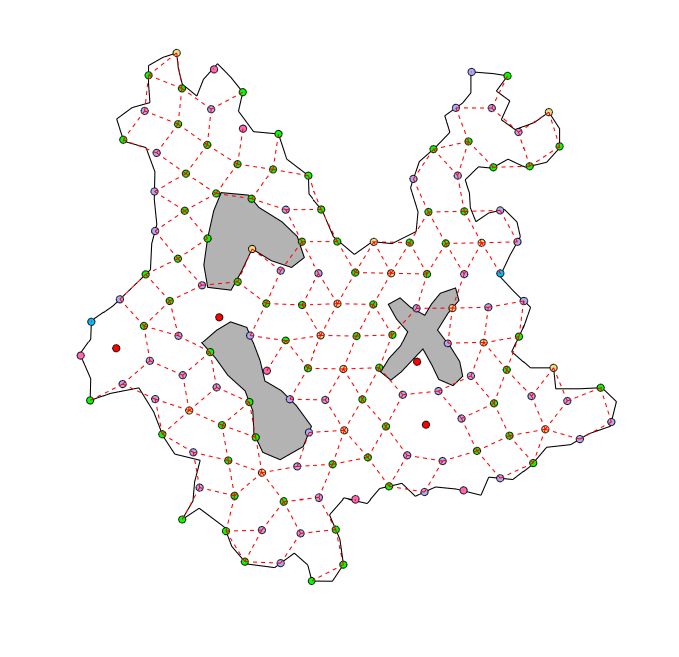}}
\subfigure[E196H5 for $p=150$]{\includegraphics[width=1.6in]{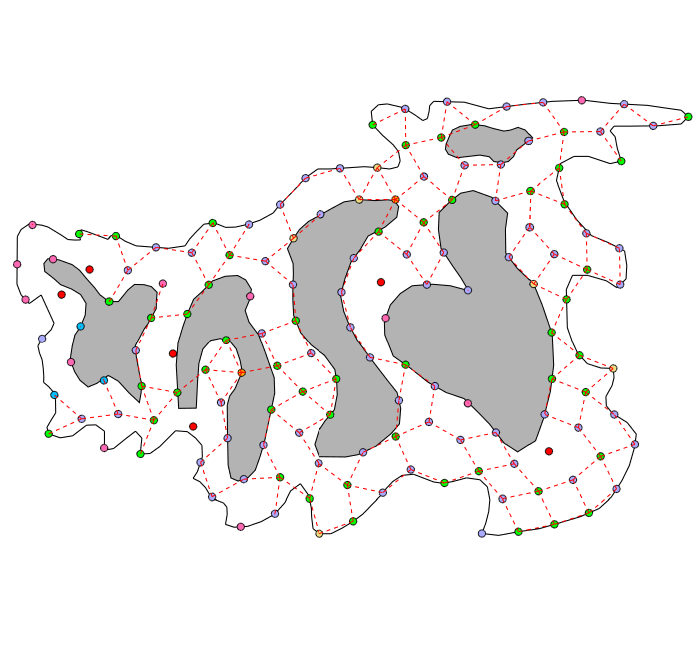}}
\subfigure[E203H3 for $p=150$]{\includegraphics[width=1.6in]{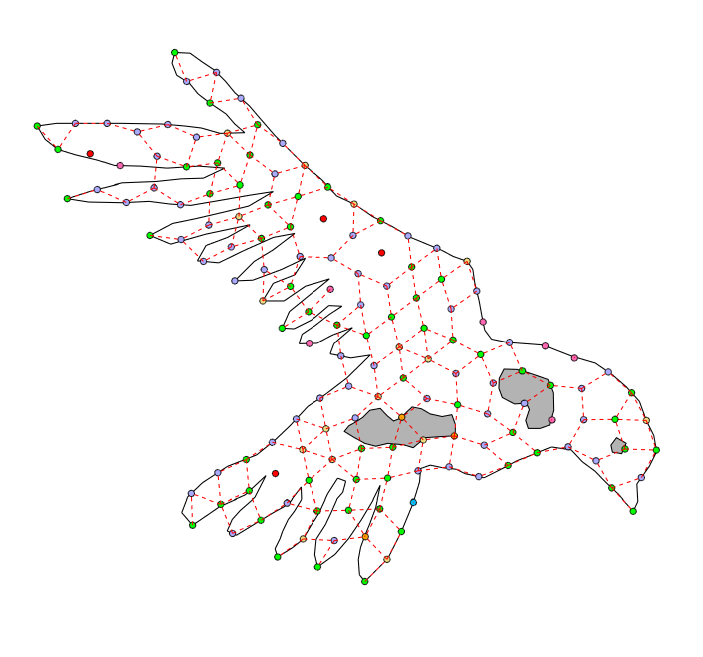}}
\caption{Best configurations found in this work for 9 representative instances, where two dispersion points are connected by a dotted line if the distance between them equals the minimum distance found $D_{best}$.}
\label{point_packing}
\end{figure*}

The objective of this section is to assess the performance of the TSGO algorithm on the CpDP problem without a boundary constraint, which corresponds to the point arrangement problem. As in Section \ref{resutls_withBC}, two variants of the TSGO algorithm, called the BH$^*$ and MBH$^*$ algorithms, were created by replacing the tabu search method respectively with the popular basin-hopping (BH) and monotonic basin-hopping (MBH) algorithms, while keeping other components unchanged. Then, the experiments were conducted based on two sets of benchmark instances with $p$ up to $p=150$, where the first set contains 30 relatively simple instances for which the region to be packed consists of a small number of edges and holes and the second set contains 45 much complicated instances for which the region to be packed consists of a large number of edges and holes. The TSGO, BH$^*$ and MBH$^*$ algorithms were independently executed 10 times for each instance, and the experimental results are summarized in Tables \ref{table_comparison2} and \ref{table_comparison2_complex} respectively for these two sets of instances, including the best objective value over 10 runs ($D_{best}$) and the average objective value ($D_{avg}$). Other statistic data in the table are the same as those in Table \ref{table_comparison1}. 

Table \ref{table_comparison2} shows that for the instances with a region consisting of a small number of edges the TSGO algorithm outperforms the BH$^*$ and MBH$^*$ algorithms. In terms of $D_{best}$, the BH$^*$, MBH$^*$ and TSGO algorithms obtained the best result among the compared algorithms for 16, 20 and 25 out of 30 instances, respectively. Nevertheless, the large \textit{p-values} mean there does not exist a significant difference between the MBH$^*$ and TSGO algorithms in terms of $D_{best}$. As for the average objective value $D_{avg}$, the BH$^*$, MBH$^*$ and TSGO algorithms respectively obtained the best result for 5, 14 and 25 instances, indicating that the TSGO algorithm outperforms the BH$^*$ and MBH$^*$ algorithms on this measure. Moreover, the small \textit{p-values} confirm that the differences between TSGO and the reference algorithms are statistically significant.

Table \ref{table_comparison2_complex} shows that for the instances with a complicated region the TSGO algorithm also outperforms the BH$^*$ and MBH$^*$ algorithms for each performance indicator considered. Moreover, for these complicated instances the superiority of the TSGO algorithm over the BH$^*$ and MBH$^*$ algorithms is more pronounced compared to the previous instances in Table \ref{table_comparison2}. Concretely, in terms of $D_{best}$, the BH$^*$, MBH$^*$ and TSGO algorithms obtained the best result for 10, 22 and 30 out of 45 instances, respectively. In terms of $D_{avg}$, the TSGO algorithm obtained the best result for 38 out of 45 instances, while the BH$^*$, MBH$^*$ algorithms  obtained the best result only for 0 and 7 instances, respectively.

To give an intuitive impression of the best configurations found, Fig. \ref{point_packing} provides the graphical representations of the best solutions for several representative instances. We observe that our TSGO algorithm is capable of effectively tackling these complicated instances for the point arrangement problem, including those with a complex container and holes as well as those consisting of a large number of edges.

\section{Analysis of Algorithmic Components}
\label{discussion}

We now analyze two important elements of our TSGO algorithm: the MBH procedure used to reinforce the intensified search of the tabu search method and the tabu search method without MBH procedure (denoted by TS-MBH). 

\subsection{Effectiveness of the MBH procedure with the tabu search method}
\label{diss_intensification}

\renewcommand{\baselinestretch}{0.95}\huge\normalsize
\begin{table}[!htp]\centering
\caption{Comparison between the TSGO algorithms with and without the MBH procedure on 40 representative instances, where the best results between the compared algorithms in terms of $R_{best}$, $R_{avg}$ and $R_{worst}$ are indicated in bold.} \label{Diss_MBH}
\begin{scriptsize}
\begin{tabular}{p{1.2cm}p{0.4cm}p{0.01cm}p{1.8cm}p{1.8cm}p{0.01cm}p{1.8cm}p{1.8cm}p{0.01cm}p{1.8cm}p{1.8cm}p{0.01cm}}
\hline
            &      & &   \multicolumn{2}{c}{$R_{best}$}    & &    \multicolumn{2}{c}{$R_{avg}$}    & &  \multicolumn{2}{c}{$R_{worst}$}  & \\
\cline{4-5} \cline{7-8} \cline{10-11}
Container   &  $p$  &  & TSGO$_1$  &  TSGO   &   &  TSGO$_1$  &  TSGO &&  TSGO$_1$ &  TSGO  &  \\
\hline
E6H0     & 200 &  & 0.0636129856          & 0.0636129856          &  & 0.0636129856 & 0.0636129856          &  & 0.0636129856 & 0.0636129856          &  \\
E9H2     & 200 &  & 0.2051661681          & \textbf{0.2051780567} &  & 0.2051047688 & \textbf{0.2051759042} &  & 0.2050264023 & \textbf{0.2051720545} &  \\
E11H0    & 200 &  & 0.2300454139          & \textbf{0.2300479293} &  & 0.2300337282 & \textbf{0.2300479293} &  & 0.2300317941 & \textbf{0.2300479293} &  \\
E12H0a   & 200 &  & 0.3303438395          & \textbf{0.3303513838} &  & 0.3303170860 & \textbf{0.3303506670} &  & 0.3302830989 & \textbf{0.3303497622} &  \\
E12H2    & 200 &  & 0.2334585328          & \textbf{0.2334870331} &  & 0.2334147356 & \textbf{0.2334801734} &  & 0.2333648795 & \textbf{0.2334692402} &  \\
E18H0    & 200 &  & 0.4190207600          & \textbf{0.4191458163} &  & 0.4189440349 & \textbf{0.4190839637} &  & 0.4188694748 & \textbf{0.4189449983} &  \\
E20H1    & 200 &  & \textbf{0.2172791291} & 0.2172791061          &  & 0.2171427598 & \textbf{0.2172676640} &  & 0.2170719917 & \textbf{0.2172289079} &  \\
E20H2    & 200 &  & 0.3446755533          & \textbf{0.3447056183} &  & 0.3446691686 & \textbf{0.3446855508} &  & 0.3446593038 & \textbf{0.3446820377} &  \\
E23H1    & 200 &  & 0.2101927860          & \textbf{0.2102037845} &  & 0.2101260082 & \textbf{0.2101661510} &  & 0.2100841375 & \textbf{0.2101431359} &  \\
E27H1    & 200 &  & 0.2595269410          & \textbf{0.2595897715} &  & 0.2594167574 & \textbf{0.2595486576} &  & 0.2591826977 & \textbf{0.2595101870} &  \\
E44H1    & 200 &  & \textbf{0.2816799600} & 0.2816788437          &  & 0.2814874593 & \textbf{0.2815233736} &  & 0.2813183912 & \textbf{0.2814098733} &  \\
E59H1    & 200 &  & 0.0111103125          & \textbf{0.0111178733} &  & 0.0110977208 & \textbf{0.0111156356} &  & 0.0110894422 & \textbf{0.0111114885} &  \\
E101H2   & 150 &  & 0.4244536975          & \textbf{0.4246301877} &  & 0.4243269607 & \textbf{0.4245517063} &  & 0.4241942738 & \textbf{0.4244364627} &  \\
E101H3   & 150 &  & 0.6022548526          & \textbf{0.6022981927} &  & 0.6019088101 & \textbf{0.6022931886} &  & 0.6016646481 & \textbf{0.6022903843} &  \\
E106H3   & 150 &  & 0.4994616507          & \textbf{0.4996758005} &  & 0.4994125371 & \textbf{0.4995933394} &  & 0.4993236995 & \textbf{0.4994616507} &  \\
E106H5   & 150 &  & 0.4440727796          & \textbf{0.4444091109} &  & 0.4439297990 & \textbf{0.4443309623} &  & 0.4438270365 & \textbf{0.4442520370} &  \\
E107H3   & 150 &  & 0.9872540160          & \textbf{0.9872881106} &  & 0.9870980239 & \textbf{0.9872064963} &  & 0.9869413117 & \textbf{0.9870158687} &  \\
E120H3   & 150 &  & 0.7306111023          & \textbf{0.7306445570} &  & 0.7303310221 & \textbf{0.7306432422} &  & 0.7298998979 & \textbf{0.7306324040} &  \\
E172H4   & 150 &  & 0.3327017919          & \textbf{0.3332646035} &  & 0.3322591879 & \textbf{0.3331663193} &  & 0.3319463556 & \textbf{0.3330925105} &  \\
E193H1   & 150 &  & 0.7408295517          & \textbf{0.7409137632} &  & 0.7399829585 & \textbf{0.7404794296} &  & 0.7390118484 & \textbf{0.7389815314} &  \\
E196H5   & 150 &  & 0.3991873083          & \textbf{0.3997206613} &  & 0.3984668794 & \textbf{0.3995365485} &  & 0.3981979101 & \textbf{0.3988107216} &  \\
E203H3   & 150 &  & 0.2498770122          & \textbf{0.2512178335} &  & 0.2445707242 & \textbf{0.2500754648} &  & 0.2372750368 & \textbf{0.2483444155} &  \\
E81H3    & 150 &  & 0.3802381769          & \textbf{0.3802766660} &  & 0.3800652384 & \textbf{0.3802496299} &  & 0.3799704229 & \textbf{0.3801959565} &  \\
E82H3    & 150 &  & 0.3490326073          & \textbf{0.3490791105} &  & 0.3488907406 & \textbf{0.3490722172} &  & 0.3488160016 & \textbf{0.3490514816} &  \\
E84H3    & 150 &  & 0.3579729123          & \textbf{0.3579871420} &  & 0.3576854002 & \textbf{0.3579385074} &  & 0.3575950893 & \textbf{0.3577198430} &  \\
E93H4    & 150 &  & 0.4946286407          & \textbf{0.4948698733} &  & 0.4944905677 & \textbf{0.4948247350} &  & 0.4943910940 & \textbf{0.4947154867} &  \\
E101H2   & 200 &  & 0.3707228206          & \textbf{0.3708266880} &  & 0.3707035113 & \textbf{0.3708085633} &  & 0.3706470633 & \textbf{0.3708000494} &  \\
E101H3   & 200 &  & 0.5270364876          & \textbf{0.5271473970} &  & 0.5265310680 & \textbf{0.5270959759} &  & 0.5261496935 & \textbf{0.5268497538} &  \\
E106H3   & 200 &  & 0.4382576604          & \textbf{0.4384167531} &  & 0.4381308018 & \textbf{0.4383871870} &  & 0.4379379803 & \textbf{0.4383279802} &  \\
E106H5   & 200 &  & 0.3909102040          & \textbf{0.3910559351} &  & 0.3906747575 & \textbf{0.3910204812} &  & 0.3904369149 & \textbf{0.3909799396} &  \\
E107H3   & 200 &  & 0.8601582219          & \textbf{0.8601603145} &  & 0.8596474782 & \textbf{0.8600083142} &  & 0.8592481578 & \textbf{0.8597126848} &  \\
E120H3   & 200 &  & 0.6348111376          & \textbf{0.6352456666} &  & 0.6346367668 & \textbf{0.6351595315} &  & 0.6343821202 & \textbf{0.6350429648} &  \\
E172H4   & 200 &  & 0.2933022486          & \textbf{0.2938264340} &  & 0.2931362207 & \textbf{0.2935365465} &  & 0.2928678971 & \textbf{0.2933340360} &  \\
E193H1   & 200 &  & 0.6493809623          & \textbf{0.6497271655} &  & 0.6492284841 & \textbf{0.6495499522} &  & 0.6491236153 & \textbf{0.6494284771} &  \\
E196H5   & 200 &  & 0.3528305296          & \textbf{0.3532967101} &  & 0.3523188998 & \textbf{0.3532569447} &  & 0.3509257201 & \textbf{0.3531768106} &  \\
E203H3   & 200 &  & 0.2225266116          & \textbf{0.2233159947} &  & 0.2216821324 & \textbf{0.2232194868} &  & 0.2206605249 & \textbf{0.2231246243} &  \\
E81H3    & 200 &  & 0.3322432698          & \textbf{0.3322827528} &  & 0.3320434488 & \textbf{0.3321969424} &  & 0.3318514747 & \textbf{0.3321376262} &  \\
E82H3    & 200 &  & 0.3043041539          & \textbf{0.3045495631} &  & 0.3041968418 & \textbf{0.3044844759} &  & 0.3041170347 & \textbf{0.3044037274} &  \\
E84H3    & 200 &  & 0.3143718343          & \textbf{0.3144135487} &  & 0.3143391699 & \textbf{0.3143854408} &  & 0.3142086130 & \textbf{0.3143467468} &  \\
E93H4    & 200 &  & 0.4321793416          & \textbf{0.4323742595} &  & 0.4317616485 & \textbf{0.4322075414} &  & 0.4314124896 & \textbf{0.4318942208} &  \\
\hline
\#Better &     &  & 2                     & 37                    &  & 0            & 39                    &  & 0            & 40                    &  \\
\#Equal  &     &  & 1                     & 1                     &  & 1            & 1                     &  & 0            & 0                     &  \\
\#Worse  &     &  & 37                    & 2                     &  & 39           & 0                     &  & 40           & 0                     &  \\
\textit{p-value}  &     &  & 6.64E-08     &                       &  & 5.26E-08     &                       &  & 7.18E-08     &                       &  \\
\hline
\end{tabular}
\end{scriptsize}
\end{table}
\renewcommand{\baselinestretch}{1.0}\large\normalsize

\begin{figure*}[ht]
\centering
\subfigure[E81H3  with $p=200$]{\includegraphics[width=2.8in]{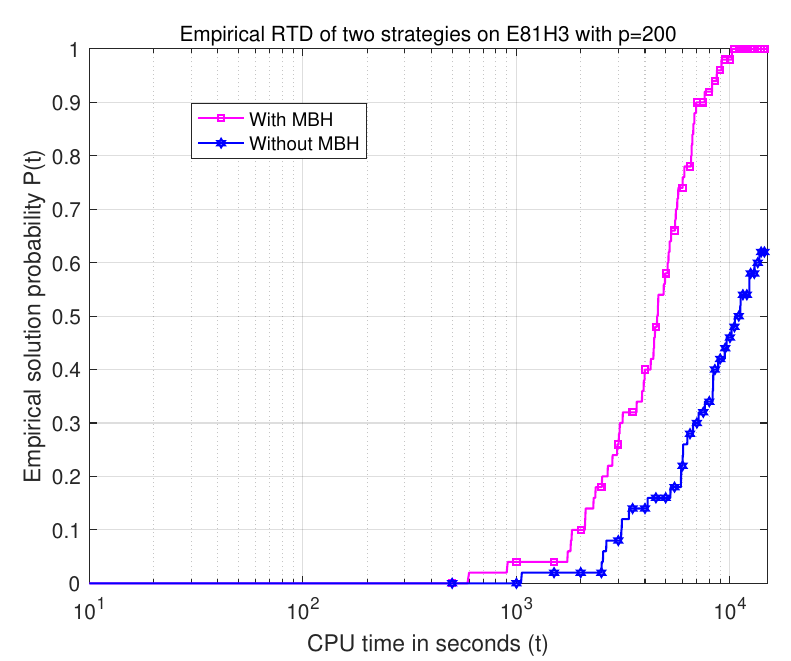}}
\subfigure[E84H3  with $p=200$ ]{\includegraphics[width=2.8in]{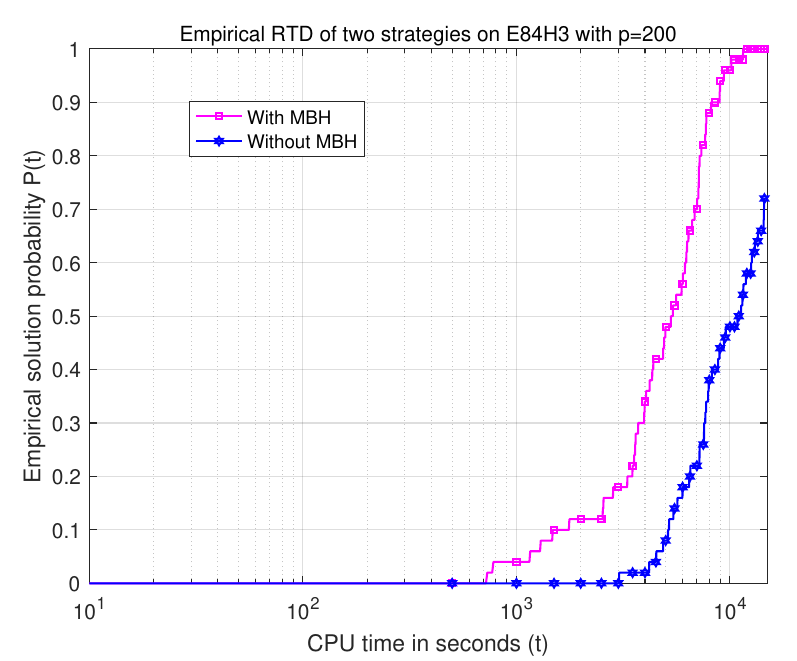}}
\subfigure[E172H4 with $p=200$]{\includegraphics[width=2.8in]{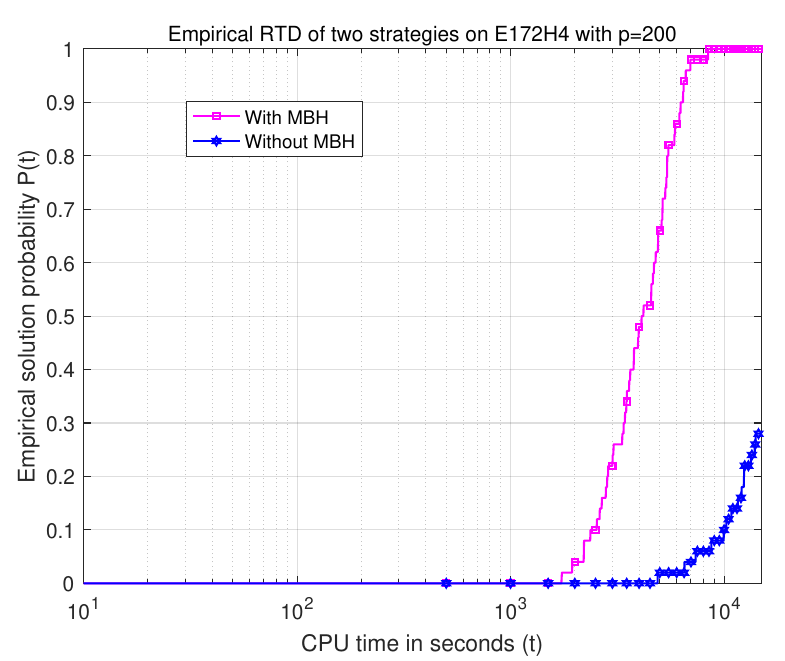}}
\subfigure[E196H5 with $p=200$]{\includegraphics[width=2.8in]{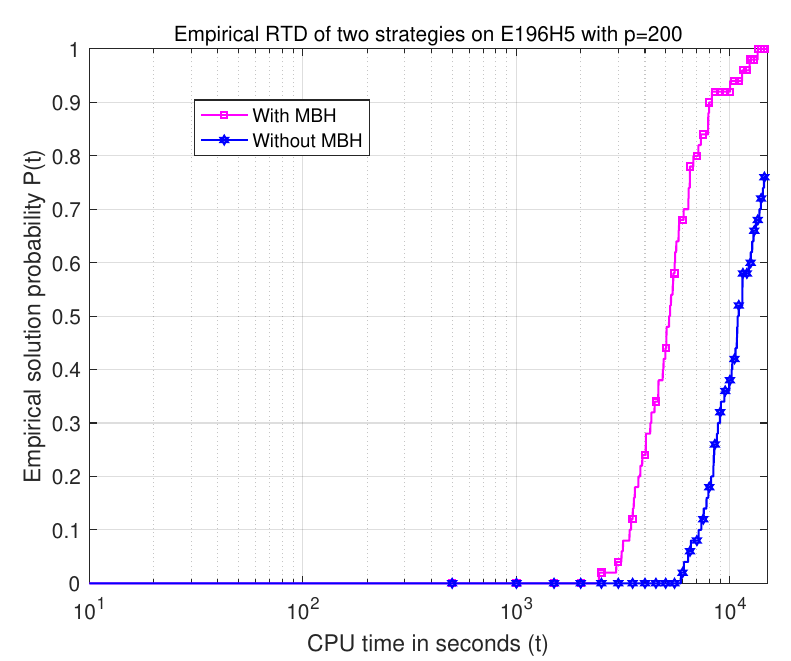}}
\caption{Empirical run-time distribution of the TSGO algorithm with and without MBH procedure on four representative instances.}
\label{analysis_MBH}
\end{figure*}

As previously noted, due to the continuity of the solution space, it is likely that a number of local minima have a similar geometric configuration that will make it difficult to distinguish among there minima using tabu search. Thus, to reinforce the intensified search, the algorithm executes a very short MBH run to improve the current solution at each iteration of the tabu search method.

To check whether the MBH procedure plays an important role for the performance of our algorithm, we carried out a comparative experiment with 40 instances of the continuous $p$-dispersion problem with boundary constraints (i.e., the equal circle packing problem). For this, we first created a variant TSGO$_1$ of our TSGO algorithm by disabling the MBH procedure while keeping other components of algorithm unchanged. Then, both of TSGO$_1$ and TSGO were run 10 times for each instance. The experimental results over 40 representative instances are summarized in Table \ref{Diss_MBH}, where columns 3-4 give the best results ($R_{best}$) over 10 independent runs respectively for TSGO$_1$ and TSGO, columns 5-6 give the average results ($R_{avg}$), and the last two columns give the worst results ($R_{worst}$). The rows ``\#Better'', ``\#Equal'', and ``\#Worse'' respectively show the numbers of instances for which the corresponding algorithm obtains a better, equal and worse result compared with the reference algorithm. The $p$-values from the Wilcoxon signed-rank test are also provided in the last row of the table.

The experimental results in Table \ref{Diss_MBH} show that the TSGO algorithm significantly outperforms its variant TSGO$_1$, which means that the MBH procedure plays an essential role in the high performance of the TSGO algorithm. Specifically, TSGO obtains a better, equal and worse result than TSGO$_1$ respectively for 37, 1 and 2 instances in terms of $R_{best}$. In terms of $R_{avg}$, TSGO performs better than TSGO$_1$ in 39 out of 40 instances and obtains equal results for the remaining instance. As for the $R_{worst}$, TSGO outperforms TSGO$_1$ on all the tested instances. Moreover, the small $p$-values ($\le 0.05$) indicate that the differences between the results of the two algorithms are statistically significant.

In addition, we employed the empirical run-time distribution (RTD) of stochastic optimization approaches to further analyze and compare the two algorithms with and without the MBH procedure. Indeed, RTD is known to be an efficient graphic representation tool to investigate the behavior of stochastic optimization algorithms \citep{hoos2000local}. For a given instance, the cumulative empirical RTD of a stochastic algorithm is a function $P(t)$ mapping the run-time $t$ to the probability of obtaining the target value within time $t$. Formally, the function $P(t)$ is defined as follows:
\begin{equation}\label{RTD_EQ}
P(t) = \frac{|\{i:rt(i)\le t\}|}{M}
\end{equation}
where $rt(i)$ represents the running time of the $i$-th successful run to obtain the target value (the average objective value $R_{avg}$ of TSGO$_1$ in Table \ref{Diss_MBH} was used in this experiment) and $M$ is the number of runs performed (where $M=50$ in this experiment).

Fig. \ref{analysis_MBH} shows TSGO's empirical RTD with and without the MBH procedure for four representative instances. One observes that for all tested instances the TSGO algorithm has a higher probability to reach the target value compared to the TSGO$_1$ algorithm without the MBH procedure under the same running time, confirming the critical role of the MBH procedure for the computational efficiency of the TSGO algorithm.

\subsection{Effectiveness of the insertion neighborhood-based tabu search method}
\label{diss_neighborhood}

\renewcommand{\baselinestretch}{1.0}\huge\normalsize
\begin{table}[!htp]\centering
\caption{Comparison between the TSGO algorithm and MBH$^*$ on 40 representative instances, where the best results between the compared algorithms in terms of $R_{best}$, $R_{avg}$ and $R_{worst}$ are indicated in bold.} \label{Diss_neighbor_moves}
\begin{scriptsize}
\begin{tabular}{p{1.2cm}p{0.4cm}p{0.01cm}p{1.8cm}p{1.8cm}p{0.01cm}p{1.8cm}p{1.8cm}p{0.01cm}p{1.8cm}p{1.8cm}p{0.01cm}}
\hline
            &      & &   \multicolumn{2}{c}{$R_{best}$}    & &    \multicolumn{2}{c}{$R_{avg}$}    & &  \multicolumn{2}{c}{$R_{worst}$}  & \\
\cline{4-5} \cline{7-8} \cline{10-11}
Container& $p$  &  & MBH$^*$               &  TSGO   &   &  MBH$^*$  &  TSGO &&  MBH$^*$ &  TSGO  &  \\
\hline
E6H0     & 200 &  & 0.0636129856          & 0.0636129856          &  & 0.0636129697          & \textbf{0.0636129856} &  & 0.0636129182          & \textbf{0.0636129856} &  \\
E9H2     & 200 &  & 0.2051540762          & \textbf{0.2051780567} &  & 0.2050474437          & \textbf{0.2051759042} &  & 0.2049556168          & \textbf{0.2051720545} &  \\
E11H0    & 200 &  & 0.2300479269          & \textbf{0.2300479293} &  & 0.2300476538          & \textbf{0.2300479293} &  & 0.2300463511          & \textbf{0.2300479293} &  \\
E12H0a   & 200 &  & 0.3303506858          & \textbf{0.3303513838} &  & 0.3303355932          & \textbf{0.3303506670} &  & 0.3302902073          & \textbf{0.3303497622} &  \\
E12H2    & 200 &  & 0.2334693785          & \textbf{0.2334870331} &  & 0.2333686152          & \textbf{0.2334801734} &  & 0.2333409856          & \textbf{0.2334692402} &  \\
E18H0    & 200 &  & 0.4191458163          & 0.4191458163          &  & 0.4189846377          & \textbf{0.4190839637} &  & 0.4188729560          & \textbf{0.4189449983} &  \\
E20H1    & 200 &  & 0.2172791061          & 0.2172791061          &  & 0.2171934541          & \textbf{0.2172676640} &  & 0.2171462736          & \textbf{0.2172289079} &  \\
E20H2    & 200 &  & 0.3446804055          & \textbf{0.3447056183} &  & 0.3446687159          & \textbf{0.3446855508} &  & 0.3446382280          & \textbf{0.3446820377} &  \\
E23H1    & 200 &  & 0.2101523773          & \textbf{0.2102037845} &  & 0.2101274162          & \textbf{0.2101661510} &  & 0.2100656909          & \textbf{0.2101431359} &  \\
E27H1    & 200 &  & 0.2595736347          & \textbf{0.2595897715} &  & 0.2594844815          & \textbf{0.2595486576} &  & 0.2592720014          & \textbf{0.2595101870} &  \\
E44H1    & 200 &  & 0.2813550521          & \textbf{0.2816788437} &  & 0.2813160515          & \textbf{0.2815233736} &  & 0.2812910866          & \textbf{0.2814098733} &  \\
E59H1    & 200 &  & 0.0111166985          & \textbf{0.0111178733} &  & 0.0111135982          & \textbf{0.0111156356} &  & 0.0111102478          & \textbf{0.0111114885} &  \\
E101H2   & 150 &  & 0.4246265745          & \textbf{0.4246301877} &  & 0.4244914723          & \textbf{0.4245517063} &  & 0.4243906129          & \textbf{0.4244364627} &  \\
E101H3   & 150 &  & 0.6022903843          & \textbf{0.6022981927} &  & 0.6022868453          & \textbf{0.6022931886} &  & 0.6022549939          & \textbf{0.6022903843} &  \\
E106H3   & 150 &  & 0.4995261676          & \textbf{0.4996758005} &  & 0.4994056398          & \textbf{0.4995933394} &  & 0.4990840235          & \textbf{0.4994616507} &  \\
E106H5   & 150 &  & 0.4443528894          & \textbf{0.4444091109} &  & 0.4442218487          & \textbf{0.4443309623} &  & 0.4440670713          & \textbf{0.4442520370} &  \\
E107H3   & 150 &  & 0.9872493423          & \textbf{0.9872881106} &  & 0.9871345218          & \textbf{0.9872064963} &  & 0.9869322205          & \textbf{0.9870158687} &  \\
E120H3   & 150 &  & 0.7306389612          & \textbf{0.7306445570} &  & 0.7304524692          & \textbf{0.7306432422} &  & 0.7300377900          & \textbf{0.7306324040} &  \\
E172H4   & 150 &  & 0.3331491441          & \textbf{0.3332646035} &  & 0.3330060576          & \textbf{0.3331663193} &  & 0.3328827019          & \textbf{0.3330925105} &  \\
E193H1   & 150 &  & 0.7407404525          & \textbf{0.7409137632} &  & \textbf{0.7406106571} & 0.7404794296          &  & \textbf{0.7402252456} & 0.7389815314          &  \\
E196H5   & 150 &  & 0.3963913171          & \textbf{0.3997206613} &  & 0.3940397475          & \textbf{0.3995365485} &  & 0.3914854411          & \textbf{0.3988107216} &  \\
E203H3   & 150 &  & 0.2474893922          & \textbf{0.2512178335} &  & 0.2423013081          & \textbf{0.2500754648} &  & 0.2402130965          & \textbf{0.2483444155} &  \\
E81H3    & 150 &  & \textbf{0.3802769198} & 0.3802766660          &  & 0.3801916372          & \textbf{0.3802496299} &  & 0.3801313844          & \textbf{0.3801959565} &  \\
E82H3    & 150 &  & 0.3490329680          & \textbf{0.3490791105} &  & 0.3489765862          & \textbf{0.3490722172} &  & 0.3488608556          & \textbf{0.3490514816} &  \\
E84H3    & 150 &  & 0.3578615502          & \textbf{0.3579871420} &  & 0.3577509179          & \textbf{0.3579385074} &  & 0.3576006015          & \textbf{0.3577198430} &  \\
E93H4    & 150 &  & \textbf{0.4948711187} & 0.4948698733          &  & 0.4948033224          & \textbf{0.4948247350} &  & 0.4946774844          & \textbf{0.4947154867} &  \\
E101H2   & 200 &  & 0.3706939791          & \textbf{0.3708266880} &  & 0.3706180722          & \textbf{0.3708085633} &  & 0.3704792462          & \textbf{0.3708000494} &  \\
E101H3   & 200 &  & 0.5270473482          & \textbf{0.5271473970} &  & 0.5265756595          & \textbf{0.5270959759} &  & 0.5263060636          & \textbf{0.5268497538} &  \\
E106H3   & 200 &  & 0.4383759370          & \textbf{0.4384167531} &  & 0.4382643170          & \textbf{0.4383871870} &  & 0.4381635633          & \textbf{0.4383279802} &  \\
E106H5   & 200 &  & 0.3909533120          & \textbf{0.3910559351} &  & 0.3908163599          & \textbf{0.3910204812} &  & 0.3907056241          & \textbf{0.3909799396} &  \\
E107H3   & 200 &  & 0.8601374785          & \textbf{0.8601603145} &  & 0.8598492621          & \textbf{0.8600083142} &  & 0.8595900817          & \textbf{0.8597126848} &  \\
E120H3   & 200 &  & 0.6351648252          & \textbf{0.6352456666} &  & 0.6350279145          & \textbf{0.6351595315} &  & 0.6349271692          & \textbf{0.6350429648} &  \\
E172H4   & 200 &  & 0.2935898258          & \textbf{0.2938264340} &  & 0.2933538090          & \textbf{0.2935365465} &  & 0.2933022486          & \textbf{0.2933340360} &  \\
E193H1   & 200 &  & \textbf{0.6497860352} & 0.6497271655          &  & 0.6494509651          & \textbf{0.6495499522} &  & 0.6493228504          & \textbf{0.6494284771} &  \\
E196H5   & 200 &  & 0.3532410620          & \textbf{0.3532967101} &  & 0.3526497635          & \textbf{0.3532569447} &  & 0.3521495511          & \textbf{0.3531768106} &  \\
E203H3   & 200 &  & 0.2217077066          & \textbf{0.2233159947} &  & 0.2213651594          & \textbf{0.2232194868} &  & 0.2209720181          & \textbf{0.2231246243} &  \\
E81H3    & 200 &  & 0.3321164693          & \textbf{0.3322827528} &  & 0.3320201998          & \textbf{0.3321969424} &  & 0.3319349657          & \textbf{0.3321376262} &  \\
E82H3    & 200 &  & 0.3044994679          & \textbf{0.3045495631} &  & 0.3044499252          & \textbf{0.3044844759} &  & 0.3043713087          & \textbf{0.3044037274} &  \\
E84H3    & 200 &  & 0.3143622086          & \textbf{0.3144135487} &  & 0.3142096904          & \textbf{0.3143854408} &  & 0.3141190415          & \textbf{0.3143467468} &  \\
E93H4    & 200 &  & 0.4319088649          & \textbf{0.4323742595} &  & 0.4317124215          & \textbf{0.4322075414} &  & 0.4315143568          & \textbf{0.4318942208} &  \\
\hline
\#Better &     &  & 3                     & 34                    &  & 1                     & 39                    &  & 1                     & 39                    &  \\
\#Equal  &     &  & 3                     & 3                     &  & 0                     & 0                     &  & 0                     & 0                     &  \\
\#Worse  &     &  & 34                    & 3                     &  & 39                    & 1                     &  & 39                    & 1                     &  \\
\textit{p-value}  &     &  & 1.14E-06     &                       &  & 1.84E-07              &                       &  & 5.34E-07              &                       &  \\
\hline
\end{tabular}
\end{scriptsize}
\end{table}
\renewcommand{\baselinestretch}{1.0}\large\normalsize

\begin{figure*}[h]
\centering
\subfigure[E81H3  with $p=200$]{\includegraphics[width=2.8in]{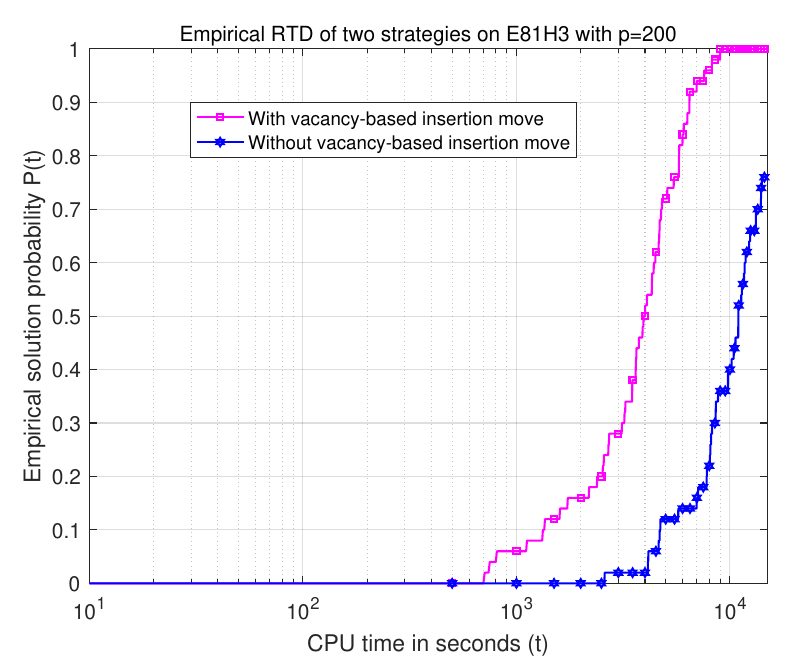}}
\subfigure[E84H3 with $p=200$ ]{\includegraphics[width=2.8in]{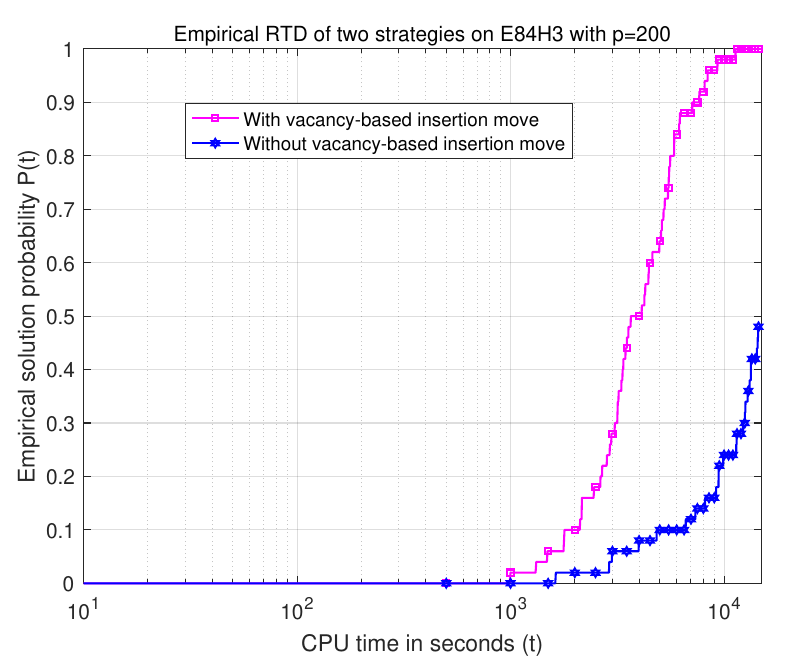}}
\subfigure[E172H4 with $p=200$]{\includegraphics[width=2.8in]{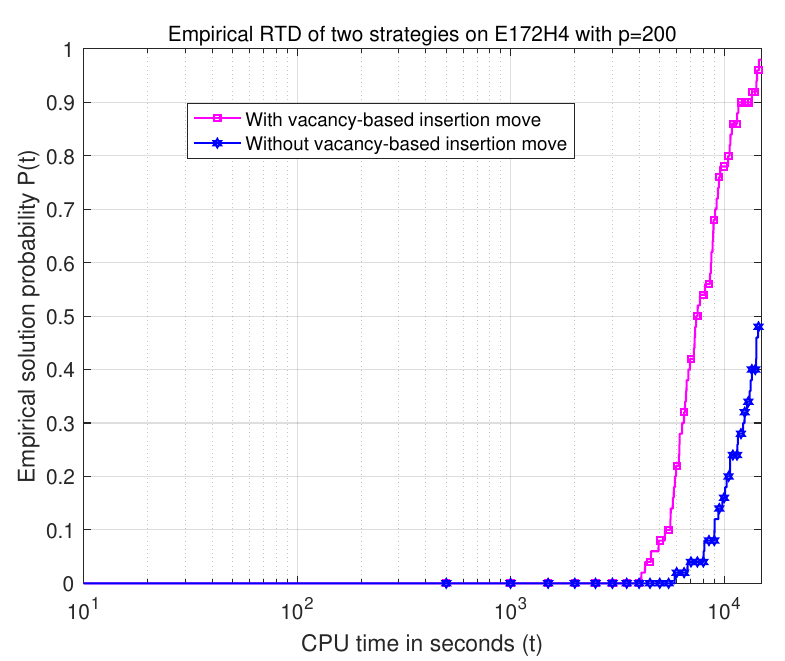}}
\subfigure[E196H5 with $p=200$]{\includegraphics[width=2.8in]{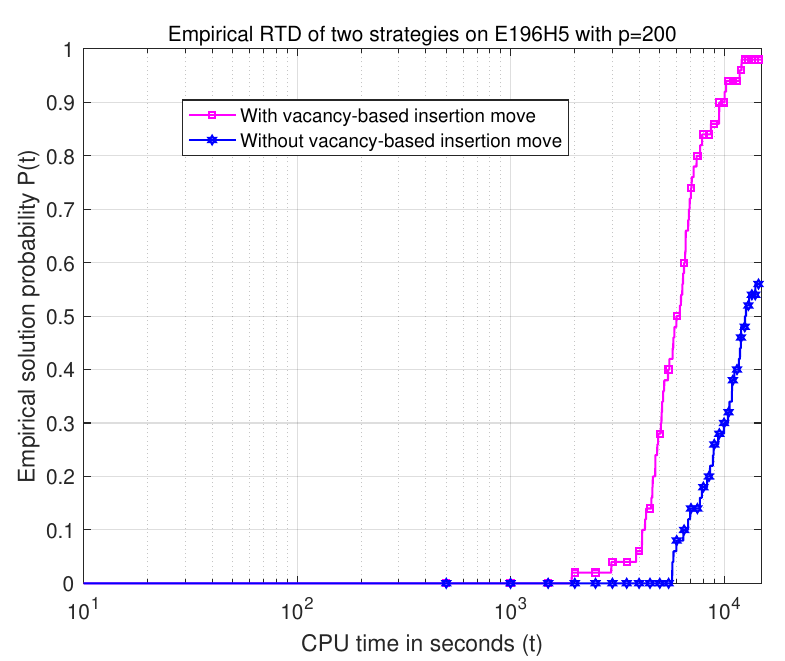}}
\caption{Empirical run-time distribution of the TSGO algorithm and MBH$^*$ on four representative instances, where MBH$^*$ corresponds to the strategy without vacancy-based insertion moves.}
\label{analysis_vacancy}
\end{figure*}

The TSGO algorithm uses the insertion neighborhood-based tabu search method without MBH (i.e., TS-MBH) as one of its main components. To check its merit, we conducted a comparative experiment based on 40 representative instances used in Section \ref{diss_intensification}. In this experiment, we created a variant of TSGO (denote by MBH$^*$) by replacing the tabu search method with the MBH method, where the search depth $\theta_{max}$ of MBH was set to its default value. Consequently, in this variant the TS-MBH method is disabled while the other components of TSGO are kept unchanged. We ran TSGO and MBH$^*$ 10 times for each instance, and the computational results are summarized in Table \ref{Diss_neighbor_moves}. We observe that TSGO outperforms MBH$^*$ in all considered performance indicators. In terms of $R_{best}$, TSGO obtained a better, equal and worse result respectively for 34, 3 and 3 instances compared to MBH$^*$. In terms of $R_{avg}$ and $R_{worst}$, TSGO obtained a better result for 39 instances and a worse result for the remaining instance. This experiment indicates that the TS-MBH method also plays a key role for the high performance of the algorithm.

To complete this comparison, Fig. \ref{analysis_vacancy} provides their empirical RTD for four representative instances, where TSGO and MBH$^*$ were performed 50 times for each instance and the target value of the RTD was set to the average objective value of MBH$^*$ in Table \ref{Diss_neighbor_moves}. We observe that TSGO has a higher probability to reach the target value compared to MBH$^*$ under the same computational time for all tested instances.

In summary, the experiments in Sections \ref{diss_intensification} and \ref{diss_neighborhood} show that disabling the MBH method or the TS-MBH method deteriorates the performance of TSGO's performance. In fact, these two component are functionally complementary, where  TS-MBH is capable of transforming the current solution into a significantly different new solution by applying the insertion moves, while MBH is capable of detecting the nearby local optimal solutions with geometrical configurations very similar to that of the current solution. Thus, the combined use of these two methods results in a high performance of the TSGO algorithm.


\section{Conclusions and Future Work}
\label{Conclusions}

The continuous $p$-dispersion problems with and without boundary constraints, which are respectively equivalent to the classic equal circle packing problem and the point arrangement problem, have numerous important real-world applications, such as the facility location and the circle cutting problems. In this study, we investigated the general cases of the continuous $p$-dispersion problems, including those that involve a non-convex multiply-connected region. For this latter general case, there does not exist an effective optimization model in the literature and it is important to address the non-trivial challenge of identifying such a model as a foundation for developing appropriate solution algorithms.

By using the penalty function method, we designed an almost everywhere differentiable optimization model for the equal circle packing problem and the point arrangement problem that compose the continuous $p$-dispersion problems with and without boundary constraints. Based on this, we proposed a global optimization method called the TSGO utilizing the proposed model. The main component of the TSGO algorithm is a tabu search method to find a feasible solution for a given minimum distance $D$ between dispersion points, coupled with a distance adjustment method to maximize the minimum distance between dispersion points while maintaining feasibility. The performance of the TSGO algorithm is assessed over a variety of existing and newly generated benchmark instances.

Thanks to the proposed optimization model, the popular continuous solvers (e.g., the quasi-Newton methods) can be applied to reach the high-precision solutions of the problems. Moreover, the proposed TSGO algorithm can reach a performance that no previous approach in the literature attains. The source code of the proposed algorithm, the set of benchmark instances used, and the best solutions found will be available online for their potential real-world applications and future algorithmic comparisons.



There are several potential directions to extend the present study in future. First, for the continuous $p$-dispersion problems without a boundary constraint, the present distance adjustment method has a slow rate of convergence due to the nature of the objective function. To improve the performance of TSGO for these problems, the distance adjustment method can be replaced by the bisection method. Second, with an appropriate modification, the present model and algorithm can be extended to handle three-dimensional problems, i.e., the continuous $p$-dispersion problems in a non-convex polyhedron containing multiple holes. Third, by changing the optimization model, the proposed TSGO algorithm can be extended to other packing problems, such as packing equal circles in a larger circle and packing equal spheres in a larger sphere.

\bibliography{PECMCR.bib}

\appendix

\end{document}